\documentclass[12pt]{article}
\usepackage{latexsym,amsmath,amssymb,amsfonts}
\usepackage{setspace}
\usepackage[a4paper,top=33mm, bottom=27mm, left=35mm, right=35mm]{geometry}
\usepackage{tikz}
\usepackage{graphicx}
\usepackage{color}

\usepackage{comment}

\usetikzlibrary{arrows,positioning} 
\usetikzlibrary{decorations.pathreplacing}
\usetikzlibrary{decorations.markings} 
\usepackage{amsthm}
\usepackage{amssymb}
\usepackage{mathrsfs}
\newtheorem{theorem}{Theorem}[section]
\newtheorem{corollary}[theorem]{Corollary}
\newtheorem{lemma}[theorem]{Lemma}

\theoremstyle{definition}

\newcommand{\wvhp}{{\rm wvhp}}
\newcommand{\BP}{{\rm BP}}

\newcommand{\inu}{\in_u}

\newcommand{\pend}{{\rm pend}}
\newcommand{\Pend}{{\rm Pend}}

\newcommand{\Cross}{{\rm Cross}}
\newcommand{\Bridge}{{\rm Bridge}}

\newcommand{\aut}{{\rm aut\,}}
\newcommand{\Aut}{{\rm Aut}}

\newcommand{\E}{{\mathbb E}}

\newcommand{\bn}{\mathbf n}

\newcommand{\bS}{\mathbf S}

\newcommand{\pr}{\mathbb P}

\newcommand{\remove}[1]{}

\newcommand{\eps}{\varepsilon}
\def\wck#1 {\underline{#1}~\marginpar{\fbox{#1} {\tiny ?}}}
\def\silent#1\par{\par}
\def\text#1{\quad\mbox{#1}\quad}

\makeatletter
\renewcommand{\@seccntformat}[1]{\@nameuse{the#1}.\quad}
\makeatother

\newcommand{\cA}{\mathcal A}
\newcommand{\cB}{\mathcal B}
\newcommand{\cC}{\mathcal C}
\newcommand{\cD}{\mathcal D}
\newcommand{\cE}{\mathcal E}
\newcommand{\cF}{\mathcal F}

\newcommand{\cH}{\mathcal H}

\newcommand{\cP}{\mathcal P}

\newcommand{\cU}{\mathcal U}

\newcommand{\tA}{\widetilde{\mathcal A}}
\newcommand{\tC}{\widetilde{\mathcal C}}
\newcommand{\tF}{\widetilde{\mathcal F}}

\newcommand{\N}{\mathbb N}

\newcommand{\Frag}{{\rm Frag}}
\newcommand{\frag}{{\rm frag}}
\newcommand{\Giant}{\rm Big}
\newcommand{\giant}{\rm big}
\newcommand{\Hb}{H^{\bullet}}

\newcommand{\Kb}{K^{\bullet}}

\newcommand{\Po}{\rm Po}

\begin{document}

\title{Pendant appearances and components in random graphs from structured classes}
\author{
Colin McDiarmid\\Department of Statistics\\ University of Oxford\\ 24 - 29 St Giles'\\ Oxford\, OX1 3LB, UK.\\cmcd@stats.ox.ac.uk}

\maketitle

\begin{abstract}
We consider random graphs sampled uniformly from a structured class of graphs, such as the class of graphs embeddable in a given surface.
We sharpen 
earlier results on pendant appearances, concerning for example numbers of leaves, and we find the asymptotic distribution of components other than the giant component, under quite general conditions. 
\end{abstract}

\section{Introduction and statement of results}
\label{sec.intro}

We consider a class $\cA$ of graphs satisfying some general structural conditions,  such as are satisfied for example by the class of graphs embeddable in a given surface.  (By a \emph{class} of graphs we mean a set of graphs which is closed under graph isomorphism.) 
For each $n \in \N$ we let $\cA_n$ be the set of graphs in $\cA$ on vertex set $[n]$, and we write $R_n \inu \cA$ to mean that the random graph $R_n$ is sampled uniformly from $\cA_n$ (implicitly assumed to be non-empty).

A connected graph $H$ has a \emph{pendant appearance} in a graph $G$ if deleting some bridge from $G$ yields a component isomorphic to $H$, see Section~\ref{subsec.pend} for a full definition.
Pendant appearances in $R_n$ have been well studied, 
with numerous applications, for example concerning symmetries and vertex degrees in $R_n$, see~\cite{msw2005,msw2006} and for example~\cite{cmcd-bham,ColinSophiaProperties}.
Here we give results on pendant appearances in $R_n$ which are both more precise and more general than earlier results, see Section~\ref{subsec.pend}; and we use these results to deduce properties of the fragment of $R_n$ (the part of the graph not in the giant component), see Section~\ref{subsec.comps}.  
In particular, we see that under quite general conditions (not implying `smoothness', see for example~\cite{cmcd-rgmc}), the asymptotic joint distribution of certain components in the fragment of $R_n$ is a Boltzmann Poisson distribution, and this may allow us to find the limiting probability that $R_n$ is connected.

After presenting our main new results, 
Section~\ref{sec.intro}
ends with a brief discussion in Section~\ref{subsec.relwork} of related work.
The plan of the rest of the paper is 
that the results on pendant appearances are proved in Section~\ref{sec.pendproofs};
the results on components are proved in Section~\ref{sec.fragproofs}; 
and Section~\ref{sec.concl} contains a few concluding remarks.
\smallskip

Recall that, given a class $\cA$ of graphs, the \emph{radius of convergence} $\rho_\cA$ or $\rho(\cA)$ of $\cA$, or of the exponential generating function (egf) of $\cA$, is given by 
\[ \rho_\cA =  \left( \limsup_{n \to \infty} \,(|\cA_n|/n!)^{1/n} \right)^{-1} ,\]
so $0 \leq \rho_\cA \leq \infty$.  We say that $\cA$ has \emph{growth constant} $\gamma_\cA$ if $0<\gamma_\cA < \infty$ and
\[ (|\cA_n|/n!)^{1/n} \to \gamma_\cA \;\; \mbox{ as } n \to \infty \, :\]
in this case clearly $\gamma_\cA = 1/\rho_\cA$.
For example the class $\cP$ of planar graphs has growth constant $\gamma_\cP \approx 27.2269$, and thus has radius of convergence $\rho_\cP =1/\gamma_\cP \approx 0.0367$, see~\cite{GN2009}.  More generally, see~\cite{CFGMN2011,cmcd-rgons}, if $S$ is any given surface then the `surface class' $\cE^S$ of graphs embeddable in $S$ has the same growth constant $\gamma_\cP$.  See for example~\cite{mt2001} for background on embeddings.

\subsection{Pendant appearances in the random graph $R_n$}
\label{subsec.pend}

Let us recall the definition of a pendant appearance. Let $H$ be a connected $h$-vertex graph with a specified root vertex $r$, forming~$\Hb$ or $H^{r}$. ($H$ will always refer to a connected graph.) 
Let $G$ be a graph and let $W$ be a set of $h$ vertices in~$G$.  If there is exactly one edge in $G$ (the \emph{link} edge) between $W$ and $V(G)\, \backslash\, W$, say the edge $wx$ with $w \in W$ and $x \in V \backslash W$, and there is an isomorphism $\phi$ from $H$ to the induced subgraph $G[W]$ such that $\phi(r)=w$, then we say that $G$ has a \emph{pendant appearance} of $\Hb$ on $W$. The set of pendant appearances of $\Hb$ in $G$ is denoted by $\Pend(G, \Hb)$, with size $\pend(G, \Hb)$.  (Thus $\pend(G, \Hb)$ is the number of sets $W$ such that $G$ has a pendant appearances of $\Hb$ on $W$.)

Let $\cA$ be a class of graphs and let $\Hb$ be a vertex-rooted connected graph.  We say that $\Hb$ is \emph{attachable to} $\cA$ if for every graph $G$ the following holds: if there is a pendant appearance of $\Hb$ in $G$ on a set $W$ and $G \,\backslash \, W \in \cA$, then $G \in \cA$.
Also, $\Hb$ is \emph{detachable from} $\cA$ if for every graph $G \in \cA$, if there is a pendant appearance of $\Hb$ in $G$ on a set $W$ then $G \, \backslash\, W \in \cA$.
For example, if $\Hb$ is a vertex-rooted connected planar graph, then $\Hb$ is attachable to and detachable from  the class $\cE^S$ of graphs embeddable in the surface $S$.
Our results are stated under weaker conditions, when $\Hb$ is `weakly attachable' or `weakly detachable', which allows the results to be applied for example in Corollary~\ref{cor.cycle} below and in~\cite{cmcd-badd-Delta} and~\cite{ColinSophiaProperties}: we postpone the corresponding definitions until just before Theorem~\ref{thm.rhoH} below.

The original Pendant Appearances Theorem~\cite{msw2005,msw2006}, see also \cite{cmcd-rgmc,cmcd-connwba2012,cmcd-rgwmc2013} says essentially that, if the class $\cA$ of graphs has a growth constant, and the vertex-rooted connected graph $\Hb$ is attachable to $\cA$, then there exists $\beta>0$ such that for $R_n \inu \cA$
\begin{equation} \label{eqn.pendapp-old}
\pr( \pend(R_n,\Hb) \leq \beta\, n) = e^{-\Omega(n)}.
\end{equation} 
A more general version of the Pendant Appearances Theorem is given and used in~\cite{cmcd-bham}:
this version assumes for $\cA$  only that $0<\rho_\cA < \infty$, and says that there exists $\beta>0$ such that
\begin{equation} \label{eqn.pendapp-old2}
  \rho\,(\{G \in \cA : \pend(G,\Hb) \leq \beta \, v(G)\,\}) > \rho_\cA 
\end{equation} 
(which easily implies~(\ref{eqn.pendapp-old}) when $\cA$ has a growth constant).  

Our main result in this section, the new
Pendant Appearances Theorem, Theorem~\ref{thm.rhoH}, extends this last result~(\ref{eqn.pendapp-old2}) in two directions: it has weaker (though perhaps less attractive) assumptions leading to new applications, and has much more precise conclusions.  Given $\Hb$, the theorem specifies the `optimal' value $\alpha$ corresponding to the constant $\beta$ in equation~(\ref{eqn.pendapp-old2}) above : part (a) of the theorem says that if  $\Hb$ is weakly attachable to $\cA$ then few graphs in $\cA_n$ have many less than $\alpha n$ pendant appearances of $\Hb$, and part (b) says that if  $\Hb$ is weakly detachable from $\cA$ then few graphs in $\cA_n$ have many more than $\alpha n$ pendant appearances of $\Hb$.

Now for the two postponed definitions. Let the class $\cA$ of graphs satisfy $0< \rho_\cA < \infty$, and let $\Hb$ be a vertex-rooted connected graph.
We say that $\Hb$ is \emph{weakly attachable} to $\cA$ if there exist $\delta>0$ and a set of graphs $\cB \supseteq \cA$
with $\rho_\cB = \rho_\cA$ such that, for all sufficiently large $n$ and all graphs $G \in \cA_n$, if $G'$ is formed from $G$ by simultaneously attaching at most $\delta n$ pendant copies of $\Hb$, then $G' \in \cB$. Observe that if $\Hb$ is attachable to $\cA$ then $\Hb$ is weakly attachable to $\cA$ (taking $\cB$ as $\cA$).  
Similarly, we say that $\Hb$ is \emph{weakly detachable} from $\cA$ if there exist $\delta>0$ and a set of graphs $\cB \supseteq \cA$ with $\rho_\cB = \rho_\cA$ such that, for all sufficiently large $n$ and all graphs $G \in \cA_n$, if $G'$ is formed from $G$ by simultaneously detaching at most $\delta n$ pendant copies of $\Hb$, then $G' \in \cB$. If $\Hb$ is detachable from $\cA$ then $\Hb$ is weakly detachable from $\cA$ (as before taking $\cB$ as $\cA$). 

For examples where a connected graph is weakly attachable to and weakly detachable from a given class of graphs but we cannot drop the word `weakly', see Corollary~\ref{cor.cycle} below and the paragraphs following it, and see also the brief comment in Section~\ref{sec.concl} on upper bounding the maximum vertex degree.
(These examples involve unrooted graphs as introduced immediately after Theorem~\ref{thm.rhoH} below.)

We denote the automorphism group of a graph $G$ by $\Aut\,G$ and let $\aut G$ be its size.
Similarly we let $\Aut\, \Hb$ be the group of automorphisms of a rooted connected graph $\Hb$ (that is, the subgroup of automorphisms of $H$ which fix the root vertex $r$), with size $\aut \Hb$.  

\begin{theorem} \label{thm.rhoH}
Let the class $\cA$ of graphs satisfy $0< \rho_\cA < \infty$.  Let $\Hb$ be a vertex-rooted connected graph, let $\alpha = \rho_{\!\cA}^{\; v(H)}/\aut\,\Hb$, and let $0<\eps<1$.
Then there exists $\nu>0$ depending on $\rho_\cA, \Hb$ and $\eps$ (and not on $\cA$ itself) such that the following holds. 
\begin{description}
\item{(a)} If $\Hb$ is weakly attachable to $\cA$ then
\[ \rho(\{G \in \cA : \pend(G,\Hb) \leq (1- \eps)\, \alpha \, v(G)\,\}) \geq \rho_\cA + \nu \,. \]
\item{(b)} If $\Hb$ is weakly detachable from $\cA$ then
\[ \rho(\{G \in \cA : \pend(G,\Hb) \geq (1+ \eps)\, \alpha \, v(G) \}) \geq \rho_\cA + \nu\, .\] 
\end{description} 
\end{theorem}

Next let us consider an unrooted version of this theorem, with the following natural corresponding definitions. Let $H$ be a connected graph. Given a graph $G$ and a set $W$ of vertices of $G$, we say that $G$ has a \emph{pendant appearance} of $H$ on $W$ if and only if $G$ has a pendant appearance on $W$ of $\Hb$ for some choice of root vertex. The set of pendant appearances of $H$ in $G$ is denoted by $\Pend(G, H)$, with size $\pend(G, H)$. Observe that, for any graph $G$, the total number $\pend(G)$ of pendant appearances (of any connected graphs) in $G$ is twice the number of bridges, so
\begin{equation}\label{eqn.pend}
\pend(G) \leq 2(v(G)-1)\,.
\end{equation} 

We say that $H$ is \emph{(weakly) attachable to} $\cA$ if $\Hb$ is (weakly) attachable to $\cA$ for each choice of root vertex; and similarly $H$ is \emph{(weakly) detachable from} $\cA$ if $\Hb$ is (weakly) detachable from $\cA$ for each choice of root vertex. 
We shall see that, by summing over the possible roots in $\Hb$, Theorem~\ref{thm.rhoH} directly yields the following  corresponding result for unrooted pendant appearances. 
\begin{corollary} \label{cor.rhoHnoroot} Let the class $\cA$ of graphs satisfy $0< \rho_\cA < \infty$. Let $H$ be an $h$-vertex connected graph, let $\alpha = h\, \rho_{\!\cA}^{\; h}/\aut H$, and let $0<\eps<1$.  Then there exists $\nu>0$ depending on $\rho_\cA, H$ and $\eps$ (and not on $\cA$ itself) such that the following holds.
\begin{description} \item{(a)} If $H$ is weakly attachable to $\cA$ then 
\[ \rho(\{G \in \cA : \pend(G,H) \leq (1- \eps)\, \alpha \, v(G)\,\}) \geq \rho_\cA + \nu \,. \] 
\item{(b)} If $H$ is weakly detachable from $\cA$ then 
\[ \rho(\{G \in \cA : \pend(G,H) \geq (1+ \eps)\, \alpha \, v(G) \}) \geq  \rho_\cA + \nu \, .\] \end{description} \end{corollary}

Given a sequence $A_n$ of events we say that $A_n$ occurs \emph{with high probability (whp)} if $\, \pr(A_n) \to 1$ as $n \to \infty$, and $A_n$ occurs \emph{with very high probability (wvhp)} if $\, \pr(A_n)=1 -e^{-\Omega(n)}$ as $n \to \infty$.  Thus for example equation~(\ref{eqn.pendapp-old}) says that $\pend(R_n,\Hb) > \beta\, n$ wvhp.
We close this subsection with a result which follows directly from Corollary~\ref{cor.rhoHnoroot} when we assume that $\cA$ has a growth constant. We shall use this result to deduce Theorem~\ref{thm.comps-badd} on components (presented in Section~\ref{subsec.comps}).
\begin{corollary} \label{cor.gcD}
Let the class $\cA$ of graphs have a growth constant, and let $R_n \inu \cA$. Let $H$ be an $h$-vertex connected graph, and let $\alpha = h\,\rho_{\!\cA}^{\,h}/\aut H$. Let $0<\eps<1$.
\begin{description}
\item{(a)} If $H$ is weakly attachable to $\cA$ then $\; \pend(R_n,H) > (1 -\eps)\, \alpha n$ \wvhp.
\item{(b)} If $H$ is weakly detachable from $\cA$ then $\; \pend(R_n,H) < (1 +\eps)\, \alpha n$ \wvhp.
\end{description}
\end{corollary}

In Corollary~\ref{cor.gcD}, it follows of course that, if $H$ is both weakly attachable to $\cA$ and weakly detachable from $\cA$ then
$\; |\pend(R_n,H)/n\: - \alpha | < \eps$ \wvhp;  and thus $\pend(R_n,H)/n \to \alpha$ in probability as $n \to \infty$.
For example, in the special case when $\cA$ is the class $\cE^S$ of graphs embeddable in a given surface $S$ and $H$ is a connected planar graph, it is known further that $\pend(R_n,H)$ is asymptotically normally distributed~\cite{GN2009,CFGMN2011}. See Section 6.7.2 of~\cite{Stufler-survey} for related results on `block-stable' classes of graphs with weights on the blocks.


\subsection{Components of the random graph $R_n$}
\label{subsec.comps}

We need some preliminary definitions, notation and results, concerning the fragment of a graph, bridge-addability, the Boltzmann Poisson distribution, and convergence in total variation.
\medskip

\noindent
\emph{The fragment of a graph}

The \emph{fragment} $\Frag(G)$ of a graph $G$ is the graph (possibly with no vertices) obtained from $G$ by deleting its largest component (breaking ties in some way), considered as an unlabelled graph.
We let $\frag(G)$ be the number of vertices in $\Frag(G)$.  We do not consider here the largest component of $G$, except to note that it has $v(G)-\frag(G)$ vertices.
Given a class $\cF$ of graphs and a graph $G$, we let $\Frag(G,\cF)$ be the unlabelled graph formed from the components of the fragment of $G$ which are in $\cF$, and let $\frag(G,\cF)$ be the number of vertices in this graph.
\medskip

\noindent
\emph{Bridge-addability}

Recall that a \emph{bridge} in a graph is an edge $e$ such that deleting $e$ yields one more component.
We say that a set $\cA$ of graphs is \emph{bridge-addable} if for every graph $G \in \cA$ with vertices $u$ and $v$ in distinct components, the graph $G+ uv$ is in $\cA$, where $G+uv$ denotes the graph obtained from $G$ by adding an edge between $u$ and $v$. (This definition is taken from~\cite{msw2006}, except that it was called  `weakly addable' there, and renamed `bridge-addable' in~\cite{cmcd-rgmc}.)  For example, the class of forests is bridge-addable, and for every surface $S$ the class $\cE^S$ of graphs embeddable in $S$ is bridge-addable.

Given a graph $G$ we let $\kappa(G)$ be the number of components of $G$. Suppose that $\cA$ is bridge-addable and $R_n \in_u \cA$. Then, by Theorem 2.2 of~\cite{msw2005}, for each $n \in \N$
\begin{equation} \label{eqn.kappa0}
\kappa(R_n) \leq_s 1 + \Po(1)
\end{equation} 
and so
\begin{equation} \label{eqn.kappa1}
\pr(R_n \mbox{ is connected}) \geq e^{-1}\,.
\end{equation}
The results~(\ref{eqn.kappa0}) and~(\ref{eqn.kappa1}) have been improved asymptotically.
Note that when $\cA$ is the class of forests, $\pr(R_n \mbox{ is connected}) = e^{-\frac12} +o(1)$ (as $n \to \infty$) \cite{Renyi59}.  By Theorem 2 of~\cite{Chapuy-Perarnau-2019} (see also~\cite{Chapuy-Perarnau-2020})
\begin{equation} \label{eqn.badd-conn12}
\pr(R_n \mbox{ is connected}) \geq e^{-\frac12} +o(1)
\end{equation}
as conjectured in~\cite{msw2006}; and indeed by Theorem 3 of~\cite{Chapuy-Perarnau-2019}, 
for any bridge-addable class $\cA$, for any $\eps>0$ 
\begin{equation} \label{eqn.badd-comp12}
\kappa(R_n) \leq_s 1+  \Po( \tfrac12 + \eps) \;\;\; \mbox{for $n$ sufficiently large}\,.
\end{equation}
%
Concerning the fragment, by inequality (7) of~\cite{cmcd-connwba2012}
\begin{equation} \label{eqn.frag}
 \E[\frag(R_n)] <2\,.
\end{equation}
(It is not known if we can bound $\E[\frag(R_n)]$ below 2.)

Finally here, given a class $\cA$ of graphs and a connected graph $H$, we say that $H$ is \emph{bridge-deletable in $\cA$}  
if for each graph $G \in \cA$, if $G$ has a bridge $e$ and at least one of the two new components in $G -e$ is a copy of $H$, then $G -e \in \cA$. This property occurs in the premises of Theorem~\ref{thm.comps-badd}.
on components of $\Frag(R_n)$. If $\cA$ is closed under deleting bridges then trivially each connected graph is bridge-deletable in $\cA$.
\medskip

\noindent
\emph{Boltzmann Poisson distribution (\cite{cmcd-rgmc})}

Given a class $\cA$ of graphs we let $\tA$ denote the set of unlabelled graphs in $\cA$, where by convention $\tA$ contains the empty graph~$\phi$.
Let $\cC$ be a non-empty class of connected graphs, and let $\cF$ be the class of graphs such that each component is in $\cC$.  Let $F(x)$ and $C(x)$ be the corresponding exponential generating functions for $\cF$ and $\cC$ respectively. Fix $\rho >0$ such that $F(\rho)= e^{C(\rho)}$ is finite; and let
\begin{equation} \label{Ceqn.mudef}
  \mu(G) = \frac{\rho^{v(G)}}{\aut G} \;\; \mbox{ for each } G \in \tF 
\end{equation}
(where $\aut \phi=1$ and so $\mu(\phi)= 1$).  We normalise these quantities to give probabilities. 
Standard manipulations, see for example~\cite{cmcd-rgmc}, show that
\begin{equation} \label{Ceqn.ugen}
  F(\rho)= \sum_{G \in \tF} \mu(G) \; \mbox{ and } \; C(\rho)= \sum_{G \in \tC} \mu(G). 
\end{equation} 
It is convenient to 
denote the corresponding \emph{Boltzmann Poisson random graph} by either $\BP (\cC,\rho)$ or $\BP (\cF,\rho)$.  It is the random unlabelled graph $R$ which takes values in the countable set $\tF$, 
such that
\begin{equation} \label{eqn.BP}
\pr(R=G) = 
 e^{-C(\rho)} \mu(G)  \;\; \mbox{ for each } G \in \tF.
\end{equation} 
Observe that $\pr(R=\phi) = 
 e^{-C(\rho)}$; and for $k \in \N$
\begin{equation}\label{eqn.orderk} \pr(v(R)=k) = \sum_{G \in \cF_k} \tfrac{\aut(G)}{k!} \, e^{-C(\rho)} \tfrac{\rho^k}{\aut(G)} =  \, e^{-C(\rho)}\, |\cF_k| \tfrac{\rho^k}{k!}\,.
\end{equation}
For each graph $G$ and each connected graph $H$, let $\kappa(G,H)$ denote the number of components of $G$ isomorphic to $H$. Then, see for example~\cite[Theorem 1.3]{cmcd-rgmc}, the random variables $\kappa(R,H)$ for $H \in \tC$ are independent, with $\kappa(R,H) \sim \Po(\mu(H))$; and $\kappa(R) \sim \Po(C(\rho))$.
\medskip

\noindent
\emph{Convergence in total variation}

Let $X_1,X_2,\ldots$ and $Y$ be random variables in some probability space.  We say that $X_n$ \emph{converges in total variation} to $Y$ if for each $\eps>0$ there exists $N$ such that for every $n \geq N$ and every measurable set $A$ we have
\[ \big| \pr(X_n \in A) - \pr(Y \in A) \big| < \eps. \]
We are interested in the case where the underlying probability space is countable, and all subsets are measurable.
\medskip

Now we present our general theorem on components, Theorem~\ref{thm.comps-badd}, followed by some applications. 
Corollary~\ref{cor.gcD} provides a key step in the proof of the theorem.
\begin{theorem} \label{thm.comps-badd}
Let the bridge-addable class $\cA$ of graphs have a growth constant. 
Let $\cC$ be a class of connected graphs such that for each graph $H$ in~$\cC$, $H$ is bridge-deletable in $\cA$, $H$ is weakly attachable to $\cA$ and $H$ is weakly detachable from $\cA$.  Then $0<C(\rho_\cA) \leq \frac12$ where $C(x)$ is the exponential generating function of $\cC$; and for $R_n \inu \cA$, 
$\, \Frag(R_n,\cC)$ converges in total variation to $\BP(\cC,\rho_\cA)$ as $n \to \infty$.
\end{theorem}
Writing $\rho$ for $\rho_\cA$, and letting $R \sim \BP(\cC,\rho)$, it follows directly from Theorem~\ref{thm.comps-badd} that
\[ \pr(\Frag(R_n,\cC)=\emptyset) \, \to \, \pr(R=\emptyset) = e^{-C(\rho)} \;\; \mbox{ as } n \to \infty. \]
Thus if also each component of $\Frag(R_n)$ is in  $\cC$ whp, then
\begin{equation} \label{eqn.Frag2}
  \pr(R_n \mbox{ is connected}) \to e^{-C(\rho)} \;\; \mbox{ as } n \to \infty\,,  
\end{equation}
see Corollaries~\ref{cor.addable},~\ref{thm.fixedS} and~\ref{cor.cycle} below.  More generally, $\kappa(\Frag(R_n,\cC))$ converges in total variation to $\kappa(R) \sim \Po(C(\rho))$; and so by~(\ref{eqn.kappa0}) each moment of $\kappa(\Frag(R_n,\cC))$ converges to the corresponding moment of $\kappa(R)$, and in particular
\[ \E[\kappa(\Frag(R_n,\cC))] \to \E[\kappa(R)] = C(\rho) \;\; \mbox{ as } n \to \infty\,. \]
We may also see easily using~(\ref{eqn.frag}) that
\begin{equation} \label{eqn.Cdash}
2 \geq \E[v(R)] = \rho \,C'(\rho)
\end{equation} 
(see immediately after the proof of Theorem~\ref{thm.comps-badd} in Section~\ref{sec.fragproofs}).

Suppose that $\cA$ is the class of forests and $\cC$ is the class of trees (both with growth constant $e$).  Then $\cA$ and $\cC$ satisfy the conditions in Theorem~\ref{thm.comps-badd}, and $\pr(R_n \mbox{ is connected}) \to e^{-\frac12}$ as $n \to \infty$ (\cite{Renyi59}, see also \cite{Moon70}): thus $C(\rho_\cA)=\frac12$, which shows that the upper bound $C(\rho_\cA) \leq \tfrac12$ in Theorem~\ref{thm.comps-badd} is best possible.
\smallskip


Recall that a \emph{minor} of a graph $G$ is a graph obtained from a subgraph of $G$ by contracting edges (and removing any loops or multiple copies of edges formed).  A class $\cA$ of graphs is \emph{minor-closed} if each minor of a graph in $\cA$ is in $\cA$. Given a minor-closed class $\cA$ of graphs, an \emph{excluded minor} for $\cA$ is a graph $G$ not in $\cA$ but such that each proper minor of $G$ is in $\cA$.  By the Robertson-Seymour Theorem, for any minor-closed class the set of excluded minors is finite (see 
for example Diestel~\cite{Diestel}).

A class of graphs is \emph{decomposable} if a graph is in the class if and only if each component is.  Following~\cite{msw2006} we call a class of graphs \emph{addable} if it is bridge-addable and decomposable. It is easy to see that a minor-closed class is addable if and only if its excluded minors are 2-connected.  
Examples of addable minor-closed classes include forests, series-parallel graphs, outerplanar graphs and planar graphs. Let $\cA$ be an addable minor-closed class of graphs.  Then $\cA$ has a growth constant~\cite{msw2006,cmcd-rgmc} (indeed perhaps this is true without assuming that the excluded minors are 2-connected, see~\cite{BNW}); and each connected graph in $\cA$ is attachable to $\cA$ and detachable from $\cA$.
  
Thus from Theorem~\ref{thm.comps-badd} and~(\ref{eqn.Frag2}) we obtain the following result, see Theorem 1.5 and Corollary 1.6 of~\cite{cmcd-rgmc}. 
\begin{corollary} \label{cor.addable}
Let $\cA$ be an addable minor-closed class of graphs, and let $R_n \inu \cA$.  Let $\cC$ be the subclass of connected graphs in $\cA$, with exponential generating function $C(x)$, and note that $\rho_\cA=\rho_\cC :=\rho$.  Then $0< C(\rho) \leq \frac12$, and $\,\Frag(R_n)$ converges in total variation to $BP(\cC, \rho)$ as $n \to \infty$; and in particular
\begin{equation} \label{eqn.addconn}
 \pr(R_n \mbox{ is connected}) \to e^{-C(\rho)} 
\;\; \mbox{ as } \; n \to \infty\,. 
\end{equation}
\end{corollary}
\noindent
When some excluded minors are not 2-connected, the behaviour of $R_n$ can be quite different, see~\cite{BM-KW-2014}.  For a different general approach, see~\cite{Stufler-Gibbs}, in particular Section 4, and see also~\cite{Stufler-survey}, in particular Section 6.7.2.
\smallskip

Consider again the class $\cP$ of planar graphs, with exponential generating function $P(x)$. Let
\[p^* = P(\rho_\cP)^{-1} 
\approx 0.96325\,, \] 
see~\cite{GN2009}.
It follows from Theorem~\ref{thm.comps-badd} or from Corollary~\ref{cor.addable} that for $R_n \inu \cP$, the fragment $\Frag(R_n)$ converges in total variation 
to $BP(\cP, \rho_\cP)$ as $n \to \infty$; and thus 
$\,\pr(R_n \mbox{ is connected}) \to p^*$ as $n \to \infty$, as shown in~\cite{GN2009}.
 
These results concerning $\cP$ extend to general surfaces.  As mentioned above (just before Section~\ref{subsec.pend}), if $S$ is any given surface then the surface class $\cE^S$ of graphs embeddable in $S$ has growth constant $\gamma_\cP$. Also, for $R_n \inu \cE^S$ whp $\Frag(R_n)$ is planar, see the start of the proof of Theorem 5.2 in~\cite{CFGMN2011}, and see Theorem 1 of~\cite{ColinSophiaProperties} for a stronger version of this result which is `uniform over all surfaces $S$'.  Thus from Theorem~\ref{thm.comps-badd} and~(\ref{eqn.Frag2}) (with $\cA=\cE^S$ and $\cC$ as the class of connected graphs in $\cP$) we obtain the next result, which essentially contains Theorems 5.2 and 5.3 of~\cite{CFGMN2011}.
\begin{corollary}
\label{thm.fixedS}
Let $\cA$ be a surface class $\cE^S$, and let $R_n \inu \cA$.  Then $\Frag(R_n)$ converges in total variation to $BP(\cP, \rho_\cP)$ as $n \to \infty$; and in particular
\begin{equation} \label{eqn.topstar}
\pr(R_n \mbox{ is connected}) \to p^* \;\; 
\mbox {as } n \to \infty\,.
\end{equation}
\end{corollary}

In results like Corollaries~\ref{cor.addable}  and~\ref{thm.fixedS}, we may be able to handle additional constraints on the random graph $R_n$; for example, an upper bound on the maximum length of a cycle (the circumference). Given a graph class $\cA$ and $t=t(n) \geq 0$,
let $\cA^t$ denote the class of graphs $G \in \cA$ such that each cycle in $G$ has length at most $t(n)$ where $n$ is the order of $G$. 
%
\begin{corollary} \label{cor.cycle}
Let $t=t(n) \to \infty$ as $n \to \infty$.
If $\cA$ is an addable minor-closed class then $\cA^t$ has growth constant $\gamma_\cA$; and for $R_n \inu \cA^t$,  $\Frag(R_n)$ converges in total variation to $BP(\cA, \rho_\cA)$ as in Corollary~\ref{cor.addable} (and so~(\ref{eqn.addconn}) holds).
If $\cA$ is a surface class $\cE^S$ then $\cA^t$ has growth constant $\gamma_\cP$; and for $R_n \inu \cA^t$, $\Frag(R_n)$ converges in total variation to $BP(\cP, \rho_\cP)$ as in Corollary~\ref{thm.fixedS} (and so~(\ref{eqn.topstar}) holds). 
\end{corollary}
%
%
Suppose temporarily that $t(n) \not\to \infty$ as $n \to \infty$.
Under the distribution $\BP(\cP,\rho_\cP)$ every cycle has strictly positive probability, so the conclusion 
that $\Frag(R_n)$ converges to $\BP(\cA,\rho_\cA)$ fails to hold when $\cA$ is $\cE^S$.  
Further, $\cA^t$ does not have growth constant $\gamma_\cA$: this follows for example from Corollary~\ref{cor.gcD}, 
see also Lemma 2.3 of~\cite{BNW}.
Similarly, if $\cA$ is an addable minor-closed class which
contains arbitrarily long cycles, then the conclusions in Corollary~\ref{cor.cycle} fail.

Let us see how to use Theorem~\ref{thm.comps-badd} to prove Corollary~\ref{cor.cycle}.  Suppose that $\cA$ is $\cE^S$ - the other case is similar except easier since then we know that $\Frag(R_n)$ is in $\cA$.
A first step is to prove that $\cA^t$ has growth constant $\gamma_\cP$, and a second step is to prove that $\Frag(R_n)$ is planar whp:  we give these proofs at the end of Section~\ref{sec.fragproofs}.
Let us assume these results for now, and quickly complete the proof. 
Observe that the class $\cA^t$ is closed under adding or deleting bridges. 
Let $H$ be a connected $h$-vertex planar graph. Taking $\cB=\cA$ and say $\delta =1$ shows that $H$ is weakly attachable to $\cA^t$;
and 
taking $\cB=\cA$ and say $\delta = 1/(2h)$ shows that $H$ is weakly detachable from $\cA^t$.
(We need the `weakly' both times here.)
Hence by Theorem~\ref{thm.comps-badd}, as $n \to \infty$ $\Frag(R_n)$ converges in total variation to $\BP(\cP, \rho_\cP)$, exactly
as in Corollary~\ref{thm.fixedS} when there is no upper bound on cycle lengths.

\smallskip

The papers~\cite{ColinSophiaSizes,ColinSophiaProperties} concern classes of graphs embeddable in order-dependent surfaces.
Given a `genus function' $g =g(n)$ we let $\cE^g$ denote the class of graphs $G$ which can be embedded in a surface of Euler genus at most $g(n)$ where $n$ is the order of $G$. Lemma 16 of~\cite{ColinSophiaProperties} shows that, when $g$ does not grow too quickly with $n$, each connected planar graph $H$ is weakly attachable to and weakly detachable from $\cE^g$.  
Thus we can use Theorem~\ref{thm.comps-badd} to extend Corollary~\ref{thm.fixedS} to this case, see Theorem~2\,(a) of~\cite{ColinSophiaProperties}.
The paper~\cite{ColinSophiaProperties} 
also contains further applications of results in the present paper.

\subsection{Related work}
\label{subsec.relwork}

Following~\cite{msw2006} (see also~\cite{msw2005}) we call a class $\cA$ of graphs \emph{smooth} if $|\cA_n|/(n |\cA_{n-1}|)$ tends to a limit as $n \to \infty$ (which must be the growth constant $\gamma_\cA$).
Bender, Canfield and Richmond~\cite{bcr08} showed in 2008 that the class $\cE^S$ of graphs embeddable in a given surface $S$ is smooth, and this also follows from the later asymptotic formula for $|\cE^S_n|$ in~\cite{CFGMN2011}.
In contrast, the graph classes considered above
(for example in Theorem~\ref{thm.comps-badd}) need not be smooth.
(In Theorem~\ref{thm.comps-badd} we could say take $\cA$ as the class of graphs $G$ which are embeddable in the sphere $\bS_0$ if $v(G)$ is even and in the torus $\bS_1$ if $v(G)$ is odd.  In this case we could also apply Theorem~\ref{thm.comps-badd} with $\cF$ as the class of planar graphs such that each component has even order.)  See also~\cite{cmcd-PAC} which is an earlier longer version of the present paper, and which contains also a version of Theorem~\ref{thm.comps-badd} with amended assumptions.

The composition method from~\cite{bcr08} can be used to show that certain other classes $\cA$ of graphs are smooth, see~\cite{cmcd-rgstrcl}; and for a class $\cA$ which is smooth, we can prove results like those presented above, in particular like Corollary~\ref{cor.gcD}, 
Theorem~\ref{thm.comps-badd}, and Corollaries~\ref{cor.addable} to~\ref{cor.cycle}. This approach naturally yields also other results concerning for example the typical size of the core, see \cite{cmcd-rgmc,cmcd-rgwmc2013,cmcd-bham,cmcd-rgstrcl}.


\section{Proofs for results on pendant appearances}
\label{sec.pendproofs}

In this section we prove 
Theorem~\ref{thm.rhoH} and then 
deduce Corollary~\ref{cor.rhoHnoroot}.  We have already noted that Corollary~\ref{cor.gcD} follows directly from Corollary~\ref{cor.rhoHnoroot}.

\begin{proof}[Proof of Theorem~\ref{thm.rhoH} (a)]
Let $\Hb$ be weakly attachable to $\cA$.  Thus there exist $\delta_0>0$ and a set $\cD$ of graphs with $\rho_\cD = \rho_\cA$ such that, for all sufficiently large $n$ and all graphs $G \in \cA_n$, if $G'$ is formed from $G$ by simultaneously attaching at most $\delta_0 n$ pendant copies of $\Hb$, then $G' \in \cD$.

Let
\[ \cB^- = \{G \in \cA : \pend(G,\Hb) \leq (1-\eps)\, \alpha \, v(G)\,\}.\]
We must show that $\rho(\cB^-) \geq \rho +\nu$ for a suitable $\nu>0$.  The idea of the proof is that from each graph in $\cA_n$ we can construct many graphs in $\cD_{n+\delta n}$ by simultaneously attaching pendant copies of $\Hb$, and if we start from graphs in $\cB^-_n$ there is limited double-counting.  Thus if $\cB^-_n$ were not small then we could form too many graphs in $\cD_{n+\delta n}$.

Let $h=v(H)$ and let $\delta =  \min \{\eps \, \alpha/4, h\,\delta_0\}$.
Let $\eta>0$ be sufficiently small that
\begin{equation} \label{eqn.eta}
  (1+\eta)^{1+\delta} (1-\eps/2)^{\delta/h} \leq 1-\eta\,,
\end{equation}
and let $\nu = \eta \rho$.
Let $k=k(n) = (\delta/h)\, n \;(\leq \delta_0 n)$, and suppose for convenience that $k \in \N$ (and so also $\delta n \in \N$).  Let $m=n+\delta n = n+hk$. Since $\cD$ has radius of convergence $\rho$, there exists $n_0$ (depending on $\cD$) 
such that 
\begin{equation} \label{eqn.rhobounds} 
\left(|\cD_n|\, / \, n!\right)^{1/n} \leq (1+\eta) \, \rho^{-1} \;\; \mbox{ for each } n \geq n_0\,.
\end{equation}

Let $n \geq n_0$.  For each $n$-subset $W$ of $[m]$ and each $n$-vertex graph $G$ in $\cB^-$ on $W$ we construct many graphs $G'$ in $\cD_{m}$ as follows.  
The first step is to partition the set $[m] \setminus W$ of the $\delta n = hk$ `extra' vertices into $k$ $h$-sets; then for each of these $h$-sets $U$ we 
put a copy of $\Hb$ on $U$ 
and add a link edge between the root of 
$\Hb$ and a vertex in $W$.
For each choice of $U$, the number of ways of doing these last steps is $\frac{h!}{\aut \Hb} \, n$.
Hence the number of constructions of graphs $G' \in \cD_m$ is
\begin{equation} \label{eqn.constrs}
\binom{m}{\delta n}\, |\cB^-_n|\, \frac{(\delta n)!}{(h!)^k k!} \left(\frac{h!}{\aut  \Hb} \right)^k n^k
= m!\, \frac{|\cB^-_n|}{n!}\, ( \aut \Hb )^{-k} \frac{n^k}{k!}.
\end{equation}

How often can each graph $G'$ be constructed?  Suppose that we start with a graph $G_0$, pick a vertex $w \in V(G_0)$, add a set $U$ of $h$ new vertices, put a copy $A$ of $\Hb$ on $U$, and add the link edge between $w$ and the root of $\Hb$, 
thus forming the new graph $G_1$. We claim that
\begin{equation} \label{claim.nd}
\pend(G_1, \Hb) \leq \pend(G_0,\Hb) + h+1\,.
\end{equation}

Let us prove this claim.
Suppose that the component of $G_0$ containing $w$ has order $r$.  If $1 \leq r \leq h-1$ then each `extra' pendant copy of $\Hb$ (other than the new one we aimed to construct) must have link edge a bridge $e$ of $\Hb$ 
such that if we delete $e$ from $\Hb$ then one of the two components formed 
has order~$r$. If $r=h$ there can be just one extra copy (with the same link edge as the constructed copy but in the opposite orientation).
If $r >h$ there can be no extra copies.  Thus if $H$ has $b$ bridges (where $0 \leq b \leq h-1$)
\[  \pend(G_1, \Hb) \leq \pend(G_0,\Hb) + \max\{2, 1+ b\} \leq  \pend(G_0,\Hb) + 
h+1\,, \]
as required in~(\ref{claim.nd}).

By~(\ref{claim.nd}), for each graph $G'$ constructed we have
\begin{eqnarray*}
\pend(G',\Hb) & \leq &
(1-\eps)\alpha n + (h+1)k\\
& \leq &
\big( (1-\eps) \alpha + 2 \delta \big)\, n
\;\;\;\; \mbox{ since } (h+1)k \leq 2hk = 2 \delta n \\
& \leq &
(1-\eps/2) \, \alpha\, n \;\;\;\; \mbox{ since } \delta \leq \eps \alpha/4.
\end{eqnarray*}
It follows that each graph $G'$ can be constructed at most
\[ \binom{\lfloor (1-\eps/2) \alpha n \rfloor}{k} \leq ((1-\eps/2) \alpha)^k \, n^k / k! \] 
times. Hence, using~(\ref{eqn.constrs}) and cancelling factors $n^k/k!$ 
\begin{eqnarray*}
\frac{|\cD_{m}|}{m!}  & \geq &
\frac{|\cB^-_{n}|}{n!}\, \big( (\aut \Hb)\, (1-\eps/2) \, \alpha \big)^{-k} \\
& = &
\frac{|\cB^-_{n}|}{n!}\, (1-\eps/2)^{-k} \, \rho^{-\delta n}
\end{eqnarray*}
since $((\aut \Hb)\, \alpha)^k = \rho^{hk} = \rho^{\delta n}$.
Therefore 
\begin{eqnarray*}
\left(\frac{|\cB^-_{n}|}{n!}\right)^{1/n} 
& \leq &
\left(\frac{|\cD_{m}|}{m!}\right)^{1/n} (1-\eps/2)^{\delta/h} \, \rho^{\delta}\\
& \leq &
((1+\eta)\rho^{-1})^{1+\delta} (1-\eps/2)^{\delta/h} \, \rho^{\delta} \;\;\; \mbox{ by } \; (\ref{eqn.rhobounds})\\
& = &
\rho^{-1} \, (1+\eta)^{1+\delta} (1-\eps/2)^{\delta/h}\\
& \leq &
(1-\eta) \rho^{-1} 
\end{eqnarray*}
by~(\ref{eqn.eta}).
It follows that
\[ \rho(\cB^-) - \rho \geq (1-\eta)^{-1} \rho - \rho > \eta \rho = \nu\,,\]
as required.
\end{proof}
\smallskip

\begin{proof}[Proof of Theorem~\ref{thm.rhoH} (b)]
The proof mirrors the proof of part (a).
Let $\Hb$ be weakly detachable from $\cA$.  Thus there exist $\delta_0>0$ and a set $\cD$ of graphs with $\rho_\cD = \rho_\cA$ such that, for all sufficiently large $n$ and all graphs $G \in \cA_n$, if $G'$ is formed from $G$ by simultaneously detaching at most $\delta_0 n$ pendant copies of $\Hb$, then $G' \in \cD$.

Let
\[ \cB^+ = \{G \in \cA : \pend(G,\Kb) \geq (1+\eps)\, \alpha \, v(G)\,\}.\]
We must show that $\rho(\cB^+) \geq \rho +\nu$ for a suitable $\nu>0$.
The idea of the proof is that from each graph in $\cB^+_n$ we can construct many graphs in $\cD_{n-\delta n}$ by deleting 
pendant copies of $\Hb$, and there is limited double-counting.  Thus if $\cB^+_n$ were not small then we could form too many graphs in $\cD_{n-\delta n}$.

Let $\delta = \min \{\eps \, \alpha/4, h\,\delta_0\}$, and assume wlog that $\delta \leq \frac12$.  Let $k=k(n) = (\delta/h)\, n \; (\leq \delta_0 n)$ and as before suppose for convenience that $k \in \N$.  Let $\eta>0$ be sufficiently small that
\begin{equation} \label{eqn.eta2}
  (1+\eps/2)^{-\delta/h} \, (1+\eta) \leq 1-\eta,
\end{equation}
and let $\nu = \eta \rho$.
We may assume that $n_0$ is sufficiently large that the upper bound~(\ref{eqn.rhobounds}) holds for this choice of $\eta$.  Let $n \geq n_0/(1-\delta)$. Since $(h+1)k \leq 2hk = 2 \delta n \leq \eps \alpha n/2$, we have
\begin{equation} \label{eqn.eta3}
 (1 + \eps)\alpha\, n - (h+1)k \geq  (1 + \eps/2)\alpha\,n.
\end{equation}

Consider a graph $G \in \cB^+_n$. 
Delete $k$ pendant appearances of $\Hb$ out of the at least $(1 + \eps)\alpha\, n$ pendant appearances of $\Hb$ in $G$, 
to form a graph $G' \in \cD$ on a set of $n-hk = n-\delta n$ vertices in $[n]$. 
Each time we delete a pendant copy of $\Hb$, by~(\ref{claim.nd}) the total number of pendant copies of $\Hb$ decreases by at most $h+1$.
Thus the number of constructions is at least
\[ |\cB^+_n| \, ((1+\eps)\alpha\, n\, - (h+1)k)^k  / k!  \geq
|\cB^+_n| \, ((1+\eps/2) \alpha\, n)^k / k! \]
by inequality~(\ref{eqn.eta3}).
The number of times $G'$ is constructed is at most the number of graphs in $\cA_n$ which we can construct by simultaneously attaching $k$ 
pendant appearances of $\Hb$ to $G'$,
which is 
\[\frac{(hk)!}{(h!)^k \, k!} \left(\frac{h! \, (1-\delta)n}{\aut \Hb}\right)^k \leq \frac{(hk)! \, n^k}{k!\, (\aut \Hb)^k} \]
(dropping the factor $(1-\delta)^k$). Hence the number of distinct graphs $G'$ constructed is at least
\[ |\cB^+_n| \, \frac{((1+\eps/2) \alpha n)^k}{k!} \,  \frac{k! \,(\aut \Hb)^k }{(hk)! \, n^k} = |\cB^+_n| \, \frac{(1+\eps/2)^k \rho^{\delta n}}{(\delta n)!} \]
since $\alpha^k (\aut \Hb)^k = \rho^{\delta n}$ and $hk=\delta n$. But this must be at most $\binom{n}{\delta n}\, |\cD_{n-\delta n}|$.  Therefore, by~(\ref{eqn.rhobounds})  since $(1-\delta)n \geq n_0$, 
\begin{eqnarray*}
|\cB^+_n| & \leq &
(\delta n)!\, (1+\eps/2)^{-k} \rho^{-\delta n} \, \binom{n}{\delta n} \, (n-\delta n)! ((1+\eta)\rho^{-1})^{n-\delta n}\\
& = &
n! \, \rho^{-n}\, (1+\eps/2)^{-k} (1+\eta)^{n-\delta n} \\
& \leq &
 n! \, \rho^{-n}\, (1+\eps/2)^{-k} (1+\eta)^{n} \\ 
& = &
n! \, \rho^{-n}\, \big( (1+\eps/2)^{-\delta/h} (1+\eta) \big)^n\\
& \leq &
n! \, \rho^{-n}\, (1-\eta)^n 
\end{eqnarray*}
where the last step uses~(\ref{eqn.eta2}).
Hence $\rho(\cB^+) -\rho \geq (1-\eta)^{-1} \rho - \rho > \nu$ as before, and this completes the proof.
\end{proof}

It remains in this section to 
deduce Corollary~\ref{cor.rhoHnoroot}.

\begin{proof}[Proof of Corollary~\ref{cor.rhoHnoroot} from Theorem~\ref{thm.rhoH}]
Consider a copy of $H$, and the corresponding action of $\Aut\, H$ on $V(H)$, with set $\cU$ of orbits.  Choose one root vertex $r_U$ from each orbit $U$, and note that for the rooted graph $H^{r_U}$ we have $\aut\, H^{r_U} = (\aut\, H) /|U|$.
For every graph~$G$
\begin{equation} \label{eqn.roots}
\pend(G,H) = \sum_{U \in \cU} \pend(G,H^{r_U})\,.
\end{equation}
Thus if
\[ \pend(G,H) \leq (1-\eps) \frac{h \rho^{h}}{\aut\,H} \, v(G)\,, \]
then for some orbit $U$
\[ \pend(G,H^{r_U}) \leq (1-\eps)\, \frac{|U|\, \rho^{h}}{\aut\,H} \, v(G) = (1-\eps)\, \frac{ \rho^{h}}{\aut\,H^{r_U}} \, v(G).\]
Hence 
\begin{eqnarray} && \rho(\{G \in \cA : \pend(G,H) \leq (1- \eps)\, \frac{h \rho^{h}}{\aut\, H} \, v(G)\,\})\nonumber \\ & \geq & \min_{r \in V(H)} \, \rho(\{G \in \cA : \pend(G,H^r) \leq (1- \eps)\, \frac{\rho^{h}}{\aut\, H^r} \, v(G)\,\})\,. \label{eqn.roots2} 
\end{eqnarray}
Similarly
\begin{eqnarray} && \rho(\{G \in \cA : \pend(G,H) \geq (1+ \eps)\, \frac{h \rho^{h}}{\aut\, H} \, v(G)\,\})\nonumber \\ & \geq & \min_{r \in V(H)} \, \rho(\{G \in \cA : \pend(G,H^r) \geq (1+ \eps)\, \frac{\rho^{h}}{\aut\, H^r} \, v(G)\,\})\,. \label{eqn.roots3} 
\end{eqnarray}

Now let us prove part (a) of Corollary~\ref{cor.rhoHnoroot}.
Suppose that $H$ is weakly attachable to $\cA$. Then each rooted graph $H^r$ is weakly attachable to $\cA$, and so by Theorem~\ref{thm.rhoH} there exists $\nu_r>0$ such that
\[ \rho(\{G \in \cA : \pend(G,H^r) \leq (1- \eps)\, \frac{\rho^{h}}{\aut\, H^r} \, v(G)\,\}) \geq \rho(\cA) + \nu_r \]
for each $r \in V(H)$.  Now let $\nu=\min_r \nu_r >0$, and the result (a) follows using~(\ref{eqn.roots2}).

We may prove part (b) in a very similar way.  Suppose that $H$ is weakly detachable from $\cA$. Then each rooted graph $H^r$ is weakly detachable from $\cA$, and so by Theorem~\ref{thm.rhoH} there exists $\nu_r>0$ such that
\[ \rho(\{G \in \cA : \pend(G,H^r) \geq (1+ \eps)\, \frac{\rho^{h}}{\aut\, H^r} \, v(G)\,\}) \geq \rho(\cA) + \nu_r \]
for each $r \in V(H)$.  As before let $\nu=\min_r \nu_r >0$, and the result (b) follows using~(\ref{eqn.roots3}).
\end{proof}


\section{Proofs for results 
on components}
\label{sec.fragproofs}

In this section we use the earlier results on pendant appearances (in particular Corollary~\ref{cor.gcD})
to prove Theorem~\ref{thm.comps-badd}.  This takes most of the section.  We also prove~(\ref{eqn.Cdash}) and give the two missing steps in the proof of Corollary~\ref{cor.cycle}.

We use six lemmas in the proof of Theorem~\ref{thm.comps-badd}. The first, Lemma~\ref{lem.kill}, is a preliminary lemma.  The second, Lemma~\ref{lem.abnew}, is the main workhorse:
it has two parts, (a) and (b), and its premises are rather general.  We use Corollary~\ref{cor.gcD} in its proof.  Lemma~\ref{lem.Hinew} has the same premises as Lemma~\ref{lem.abnew} when we combine the two parts of that lemma.  Indeed, each one of Lemmas~\ref{lem.Hinew}, \ref{lem.Hi2new}, \ref{lem.Crhonew} and~\ref{lem.prefree} has these same premises. Lemma~\ref{lem.prefree} and the inequalities following its proof complete the proof of Theorem~\ref{thm.comps-badd}.

Given a graph $G$ and a class $\cC$ of connected graphs, let $\Pend(G,\cC)$ denote the union of the (disjoint) sets $\Pend(G,H)$ for $H \in \cC$.  
\begin{lemma} \label{lem.kill}
Let $\cC$ be a class of connected graphs with $\max\{v(H): H \in \cC \} = h^* \in \N$.  Let the graph $G$ have a pendant appearance of a graph $H_0 \in \cC$ on the set $W \subseteq V(G)$ (so $H_0$ is the induced subgraph $G[W]$), and let $G'$ be $G \setminus W$.  Let $Q= \Pend(G,\cC) \setminus \Pend(G',\cC)$. Then $|Q| \leq 2h^*$.
\end{lemma}
\begin{proof}
Each member of $Q$ is specified by its oriented link in $G$.  Let $wx$ (with $w \in W$) be the oriented link of the pendant appearance of $H_0$ on $W$ (pointing away from $W$), and let $H_0^+$ be the graph obtained by adding the (unoriented) edge $\{w,x\}$ to $H_0$.

Let $Q'$ be the set of members of $Q$ such that the (unoriented) bridge corresponding to its link is not in $H_0^+$.  Observe that if a member of $Q'$ has vertex set $W'$ then we must have $W \cup \{x\} \subseteq W'$ and so its oriented link must point away from $W$. Consider a member of $Q'$ with vertex set $W_1$ as large as possible, and suppose that it is a pendant  appearance of $H_1 \in \cC$.  Then for every other member of $Q'$ the bridge corresponding to its link is a bridge in $H_1$ and not in $H_0^+$.  There are at most $v(H_1) - v(H_0)-1 \leq h^* -v(H_0)-1$ such bridges, so $|Q'| \leq h^* -v(H_0)$.
But there are at most $v(H_0)$ bridges in $H_0^+$, so
\[ |Q| \leq |Q'| + 2v(H_0) \leq h^*+v(H_0) \leq 2h^*,\]
as required.
\end{proof}

The bound in Lemma~\ref{lem.kill} is tight.
For example, suppose that $\cC$ consists of all paths with at most $h^*$ vertices. Let the graph $G$ have a component which is a path $P$ with $h^*+1$ vertices, and let $W$ be the set of vertices in $P$ apart from one end-vertex.  Then for each of the $h^*$ bridges in $P$, both of the orientations yield a pendant appearance of a graph in $\cC$, and thus $|Q|= 2h^*$.

Recall that $\kappa(G,H)$ is the number of components of $G$ which are isomorphic to $H$.
\begin{lemma} \label{lem.abnew}
Let the bridge-addable class $\cA$ of graphs have growth constant~$\gamma$, and let $\rho =\rho_\cA\, (=1/\gamma)$.  Let $\cC$ be a non-empty class of connected graphs, and let each graph in $\cC$ be bridge-deletable in $\cA$.
List the (unlabelled) graphs in $\tC$ as $H_1,H_2,\ldots$ so that $v(H_1) \leq v(H_2) \leq \cdots$.
Call $k \in \N$ \emph{relevant} when $k \leq |\tC|$ (where $\tC$ may be infinite).  For each relevant $i$  let $h_i=v(H_i)$, let $\alpha_i= \rho^{h_i}/ \aut H_i$ and let $X_i \sim \Po(\alpha_i)$; and let $X_1, X_2,\ldots$ be independent.  Also, for each relevant $k$ let $\sigma_k = \sum_{i=1}^k \alpha_i$. Let $R_n \inu \cA$.

(a) 
If each graph in $\cC$ is weakly attachable to $\cA$, then for every $\eps>0$ there exists a relevant $k_0 \in \N$ such that for every relevant $k \geq k_0$ the following two statements hold.\\
(i) For all $n$ sufficiently large
\begin{equation} \label{eqn.probsa0} \pr\left( \land_{i=1}^k \,\kappa(R_n,H_i) = 0 \right) \leq (1+\eps)\, e^{-\sigma_k}\,.
\end{equation}
(ii) For every $n^* \in \N$ there exists $N^* \in \N$ such that for all $\,0 \leq n_1,\ldots,n_k \leq n^*$ and all $n \geq N^*$
\begin{equation} \label{eqn.probsa}
\pr\left( \land_{i=1}^k \,\kappa(R_n,H_i) = n_i \right) \geq (1\!-\!\eps)\, 
\pr\left( \land_{i=1}^k \,\kappa(R_n,H_i) = 0 \right) e^{\sigma_k} \prod_{i =1}^{k} \pr(X_i=n_i) \,.
\end{equation}

(b) 
If each graph in $\cC$ is weakly detachable from $\cA$,
then for every $\eps>0$ and every relevant $k \in \N$ the following two statements hold.\\
(i) For all $n$ sufficiently large
\begin{equation} \label{eqn.probsb0}
\pr\left( \land_{i=1}^k \,\kappa(R_n,H_i) = 0 \right) \geq (1-\eps)\, e^{-\sigma_k}\,.
\end{equation}
(ii) For every $n^* \in \N$ there exists $N^* \in \N$ such that for all $\,0 \leq n_1,\ldots,n_k \leq n^*$ and all $n \geq N^*$
\begin{equation} \label{eqn.probsb}
\pr\left( \land_{i=1}^k \,\kappa(R_n,H_i) = n_i \right) \leq (1\!+\!\eps)\, 
\pr\left( \land_{i=1}^k \,\kappa(R_n,H_i) = 0 \right) e^{\sigma_k} \prod_{i =1}^{k} \pr(X_i=n_i) \,.
\end{equation}
\end{lemma}
\noindent
Note the asymmetry here: we require $k$ to be large in part (a) but not in part~(b).
\begin{proof}
For all $n \geq 1$, relevant $k$ and $n_1,\ldots,n_k \geq 0$ let
\[ \cA_n(n_1,\ldots,n_k) =  \{ G \in \cA_n : \land_{i=1}^k \, \kappa(G,H_i) = n_i \}. \]
Also, for relevant $k$ write $\tC_{(k)}$ for $\{H_1,\ldots,H_k\}$.  Let $0<\eps <1$.
\medskip

Proof of part (a). 
After some preliminaries, we see how to construct many graphs in $\cA_n(n_1,\ldots,n_k)$, by choosing a `typical' graph $G$ in $\cA_n({\bf 0})$ and deleting the link edge in $G$ from $n_i$ pendant appearances of $H_i$ for each $i \in [k]$.  We avoid creating unwanted extra components $H_i$ by performing these deletions only within the giant component of $G$. This procedure yields~(\ref{eqn.probsa}), and then~(\ref{eqn.probsa}) easily gives~(\ref{eqn.probsa0}). 

Assume that each graph in $\cC$ is weakly attachable to $\cA$.  Pick a relevant $k_0 \in \N$ sufficiently large that either $v(H_{k_0}) \geq 24/\eps$, or $\tC$ is finite and $k_0=|\tC|$.  Let $k \geq k_0$ be relevant, and let $h^*=v(H_k)$.  Let $n^* \in \N$. (We shall define a suitable $N^*$ below.) Let $s^*= n^* \, \sum_{i=1}^{k} h_i$. Let $\cB$ be the (bad) class of graphs $G \in \cA$ such that $\frag(G) \geq 24/\eps$.  Since  $\E[\frag(R_n)]<2$ by~(\ref{eqn.frag}) (and Markov's inequality),
\begin{equation} \label{eqn.Bprob2}
\pr(R_n \in \cB) < \eps/12 \;\; \mbox{ for each } n \in \N\,. 
\end{equation}
Note that if $G \in \cA \setminus \cB$ then by~(\ref{eqn.pend}) the total number of pendant appearances in $\Frag(G)$ is less than $48/\eps$.
Also, by~(\ref{eqn.kappa1})
\begin{equation} \label{eqn.An0}
  |\cA_n({\bf 0})| \geq |\cA_n| /e > |\cA_n| /3 \;\;\mbox{ for each } n \in \N 
\end{equation}
where ${\bf 0}$ is the $k$-tuple of 0's.

Let $t^*= k n^*$.  Let $0<\eta< \frac12$ satisfy $(1- 2 \eta)^{t^*} \geq 1- \eps/2$.
Let $\cA^-$ be the set of graphs $G \in \cA$ such that
$\pend(G,H_i) \leq (1-\eta)\, \alpha_i \, h_i\, v(G)\,$ for some $i \in [k]$.  By Corollary~\ref{cor.gcD} (a) we have $R_n \not\in \cA^-$ whp (indeed wvhp).
Let $N_1$ be sufficiently large that for all $n \geq N_1$ we have $|\cA_n^-| \leq (\eps/12) \, |\cA_n|$.
Let $N_2 = (2h^*t^* + 48/\eps)/(\eta \min_{1 \leq i \leq k} \alpha_i h_i)$, let $N_3 = 24/\eps + t^* + h^* +1$, and let $N^*$ be the maximum of $N_1$, $N_2$ and $N_3$.
\smallskip

We have now completed the preliminaries for the proof of inequality~(\ref{eqn.probsa}).  Let $\,0 \leq n_1,\ldots,n_k \leq n^*$, and let $t= \sum_{i=1}^k n_i$ (so $0 \leq t \leq t^*$).  Let $n \geq N^*$. Let
\[G \in \cA_n({\bf 0}) \setminus (\cA^- \cup \cB).\]
We shall ignore $\Frag(G)$. We aim to choose a pendant appearance in $\Giant(G)$ (the `giant' component of $G$) of a graph $H_j$ in $\tC_{(k)}$ such that $n_j \geq 1$, delete the link edge from $G$ to form the graph $G' \in \cA_n$ (note that $G' \in \cA_n$ since $H_j$ is bridge-deletable in $\cA$, and $G'$ has a component $H_j$), and update the numbers $n_i$ to $n'_i$ (setting $n'_j = n_j -1$ and $n'_i =n_i$ for $i \neq j$); then choose a pendant appearance in $G'$ less $\Frag(G)$ of a graph $H_{j'}$ in $\cC_{(k)}$ with $n'_{j'} \geq 1$, delete the link edge and so on: we continue until either we have chosen $n_i$  pendant copies of $H_i$ for each $i \in [k]$ and deleted their link edges, or we fail at some stage.  

By Lemma~\ref{lem.kill}, deleting the link edge of a pendant appearance of a graph in $\tC_{(k)}$ can destroy at most $2h^*$ pendant appearances of a graph in $\tC_{(k)}$, so in total we can destroy at most $2h^*t \leq 2h^*t^*$ pendant appearances of graphs in $\tC_{(k)}$.  Observe that $(1-\eta) \alpha_i h_i n - 2h^* t^* -48/\eps \geq  (1- 2\eta) \alpha_i h_i n$ (since $n \geq N^* \geq N_2$). Thus starting from $G$ we do not fail, and we may construct at least 
\[ \prod_{i=1}^k \frac{\big( (1-2\eta) \alpha_i h_i n \big)^{n_i}}{n_i !} = (1-2\eta)^{t} \, n^t \,\prod_{i=1}^k \frac{(\alpha_i h_i)^{n_i}}{n_i!} \]
graphs $G' \in \cA_n$.  Observe that in the whole process we delete $t \leq t^*$ vertices from $\Giant(G)$, leaving a component with at least $n-\frag(G)-t^* > h^*$ vertices (since $n \geq N^* \geq N_3$): thus there are no `extra' components in $\cC_{(k)}$, and $G' \in \cA_n(n_1,\ldots,n_k)$.

Also, each graph $G'$ is constructed at most $n^t \,\prod_{i=1}^k h_i^{n_i}$ times, since this quantity bounds the number of ways of re-attaching the components $H_i$ to form pendant appearances.  Hence
\[ |\cA_n({\bf 0}) \setminus (\cA^- \cup \cB)| \cdot (1-2\eta)^{t} \, n^t \,\prod_{i=1}^k \frac{(\alpha_i h_i)^{n_i}}{n_i!} \leq | \cA_n(n_1,\ldots,n_k)|\cdot n^t \,\prod_{i=1}^k h_i^{n_i}.\]
Thus
\[ | \cA_n(n_1,\ldots,n_k)|  \geq |\cA_n({\bf 0}) \setminus (\cA_n^- \cup \cB_n)| \; (1-2\eta)^{t} \,\prod_{i=1}^k \frac{\alpha_i^{n_i}}{n_i!}.\]
But $|\cA_n^-| \leq (\eps/12) \, |\cA_n|$ since $n \geq N^* \geq N_1$; and $|\cB_n| \leq (\eps/12) \, |\cA_n|$ by (\ref{eqn.Bprob2}); so
$|\cA_n^- \cup \cB_n| \leq (\eps/6)\, |\cA_n| \leq (\eps/2)\, |\cA_n({\bf 0})|$ by~(\ref{eqn.An0}). Thus
\begin{eqnarray*}
  | \cA_n(n_1,\ldots,n_k)| 
& \geq &
  (1-\eps/2) \, |\cA_n({\bf 0})| \; (1-2\eta)^{t} \,\prod_{i=1}^k \frac{\alpha_i^{n_i}}{n_i!}\\
& \geq &
  (1-\eps)\, |\cA_n({\bf 0})|\, e^{\sigma_k} \,\prod_{i=1}^k \pr(X_i=n_i)  
\end{eqnarray*}
since $(1-2\eta)^t \geq (1-2\eta)^{t^*} \geq 1-\eps/2$, and $(1-\eps/2)^2 \geq 1-\eps$. Dividing both sides of this inequality by $|\cA_n|$ yields inequality~(\ref{eqn.probsa}).

We next prove~(\ref{eqn.probsa0}). We may assume that $0<\eps \leq \frac13$, so $(1-\eps/3)^2 \geq (1+\eps)^{-1}$.  We use~(\ref{eqn.probsa}) with $\eps$ replaced by $\eps/3$ to prove~(\ref{eqn.probsa0}). 
Let $n^*$ be sufficiently large that $\prod_{i=1}^k \pr(X_i \leq n^*)  \geq 1-\eps/3$, let $N^*$ be as given for~(\ref{eqn.probsa}), and let $n \geq N^*$.  By summing~(\ref{eqn.probsa}) (with $\eps$ replaced by $\eps/3$) over all $0 \leq n_1,\ldots,n_k \leq n^*$ we obtain
\begin{eqnarray*}
  |\cA_n|
& \geq &
  |\{G \in \cA_n : \land_{i=1}^{k} \kappa(G,H_i) \leq n^*\}|\\ 
& \geq &
  (1-\eps/3)\, |\cA_n({\bf 0})|\, e^{\sigma_k} \,\prod_{i=1}^k \pr(X_i \leq n^*)\\
& \geq &
  (1-\eps/3)^2\, |\cA_n({\bf 0})|\, e^{\sigma_k}\,.
\end{eqnarray*}
Hence
$|\cA_n({\bf 0})| \leq (1+\eps)\, |\cA_n|\, e^{-\sigma_k}$, which yields~(\ref{eqn.probsa0}) and thus completes the proof of part~(a).
\medskip

Proof of part (b).  We may think of the inequality~(\ref{eqn.probsb}) which we need to prove as giving a lower bound on $|\cA_n({\bf 0})|$.
After some preliminaries, we see how to construct enough graphs in $\cA_n({\bf 0})$ by picking a `typical' graph $G \in \cA_n(n_1,\ldots,n_k)$ and adding a link edge between each component $H_i$ and $\Giant(G)$ (to form a pendant appearance of $H_i$).   By linking just to $\Giant(G)$ we ensure that we avoid creating any extra unwanted components in $\cC$. Finally we use~(\ref{eqn.probsb}) to prove~(\ref{eqn.probsb0}).

Assume that each graph in $\cC$ is weakly detachable from $\cA$. Let $k \in \N$ be relevant, let $h^*=v(H_k)$, and let $n^* \in \N$. Let 
\[ \beta = \min_{0 \leq n_i \leq n^*} \prod_{i=1}^k \frac{\alpha_i^{n_i}}{n_i!} \] 
where the minimum is over all $0 \leq n_1,\ldots,n_k \leq n^*$, and note that $\beta>0$.  As before, let $s^* = n^* \sum_{i=1}^k h_i$ and $t^*=k \, n^*$. Let $\eta>0$ satisfy $(1+2\eta)^{2t^*} + \eta \leq 1+\eps$.

Let $f_0 = 12/\eta \beta$, and let $\cB$ be the (bad) class of graphs $G \in \cA$ such that $\frag(G) \geq f_0$.  Since  $\E[\frag(R_n)]<2$ by~(\ref{eqn.frag}),
\begin{equation} \label{eqn.Bprob3}
\pr(R_n \in \cB) < \tfrac16 \eta \beta \;\; \mbox{ for each } n \in \N\,. 
\end{equation}

Let $\cA^+$ be the set of graphs $G \in \cA$ such that
$\pend(G,H_i) \geq (1+\eta)\, \alpha_i \, h_i\, v(G)\,$ for some $i \in [k]$.  By Corollary~\ref{cor.gcD} (b) we have $R_n \not\in \cA^+$ whp (indeed wvhp). Let $N_1$ be sufficiently large that for all $n \geq N_1$ we have
$ |\cA^+_n| \leq \frac16 \,\eta \, \beta \, |\cA_n|$;
let $N_2 = 2 t^*h^*/(\eta\, \min_{1 \leq i \leq k} \alpha_i h_i)$; let $N_3= f_0 (1+2\eta)/2\eta$; let $N_4 = f_0 + h^*$; and finally let $N^*$ be the maximum of $N_1, N_2, N_3, N_4$. This completes the preliminaries.

Let $0 \leq n_1,\ldots,n_k \leq n^*$ be given, and let $n \geq N^*$. We want to prove that inequality~(\ref{eqn.probsb}) holds.
Let $t= \sum_{i=1}^k n_i$, so $0 \leq t \leq t^*$; and let $s= \sum_{i=1}^k h_i n_i$, so $0 \leq s \leq s^*$. Let
\[ G \in \cA_n(n_1,\ldots,n_k) \setminus (\cA^+_n \cup \cB_n).\]
In $G$, for each component isomorphic to one of the graphs $H_i$ in $\cC_{(k)}$ (if $t \geq 1$) we pick one of the $h_i$ vertices in the component and one of the $n - \frag(G) \geq n- f_0$ vertices in $\Giant(G)$ and add a bridge between these two vertices, which forms a pendant copy of $H_i$ in the new largest component, of order at least $n-f_0 +1 > h^*$, which is too big to be in $\cC_{(k)}$.  Note that the new graph $G'$ is in $\cA_n$ since $\cA$ is bridge-addable.  In this way we construct at least
\[ \prod_{i=1}^k (h_i (n- f_0))^{n_i} = (n-f_0)^t \prod_{i=1}^k h_i^{n_i}\]
graphs $G' \in \cA_n({\bf 0})$. By Lemma~\ref{lem.kill}, each time we add a bridge as here we can construct at most $2h^*$ new pendant copies of graphs in $\cC_{(k)}$, so in total we construct at most $2 t h^*$  new pendant copies of graphs in $\cC_{(k)}$.  Since $n \geq N_2$ we have $2 t h^* \leq \eta\, \alpha_i h_i n$ for each~$i$, and so
\[ \pend(G',H_i) \leq (1+\eta) \alpha_i h_i n + 2 t  h^* \leq (1+2 \eta) \alpha_i h_i n ; \]
so the number of times $G'$ is constructed is at most 
\[ \prod_{i=1}^k \binom{\lfloor (1+2\eta) \alpha_i h_i n \rfloor}{n_i}
\leq \prod_{i=1}^k \frac{((1+2\eta) \alpha_i h_i n)^{n_i}}{n_i!}
= (1+2\eta)^{t} \, n^t \,\prod_{i=1}^k \frac{(\alpha_i h_i)^{n_i}}{n_i!} \,.\]
Hence 
\[ | \cA_n(n_1,\ldots,n_k) \setminus (\cA^+ \cup \cB) |\: \big(1-\tfrac{f_0}{n}\big)^t n^t \prod_{i=1}^k h_i^{n_i} \leq |\cA_n({\bf 0})| \, (1+2\eta)^{t} \, n^t \,\prod_{i=1}^k \frac{(\alpha_i h_i)^{n_i}}{n_i!}.\]
Thus, using $1-\frac{f_0}{n} \geq (1+2\eta)^{-1}$
(which holds since  $\,n \geq N^* \geq f_0 \,(1+2\eta)/2\eta\,$) we have
\[ | \cA_n(n_1,\ldots,n_k) \setminus (\cA_n^+ \cup \cB_n) |  \leq |\cA_n({\bf 0})| \; (1+2\eta)^{2t} \,\prod_{i=1}^k \frac{\alpha_i^{n_i}}{n_i!}.\]
Now, using~(\ref{eqn.Bprob3}) and then~(\ref{eqn.An0}),
\[ |\cA_n^+ \cup \cB_n| \leq \tfrac13\, \eta \beta \, |\cA_n| \leq \eta \beta\, |\cA_n({\bf 0})|\,. \]
Hence
\begin{eqnarray*}
  | \cA_n(n_1,\ldots,n_k)|  
& \leq &
  |\cA_n({\bf 0})| \; \big[(1+2\eta)^{2t} \,\prod_{i=1}^k \frac{\alpha_i^{n_i}}{n_i!} + \eta \beta \big]\\
& \leq &
  |\cA_n({\bf 0})| \; ( (1+2\eta)^{2t} + \eta) \,\prod_{i=1}^k \frac{\alpha_i^{n_i}}{n_i!}\\
& \leq &
   (1+\eps)\, |\cA_n({\bf 0})|\, e^{\sigma_k} \,\prod_{i=1}^k \pr(X_i =n_i)
\end{eqnarray*}
by our choice of $\eta$.  Thus inequality~(\ref{eqn.probsb}) holds.

Finally, we use~(\ref{eqn.probsb}) with $\eps$ replaced by $\eps/3$ to prove~(\ref{eqn.probsb0}) (following the pattern of the proof of part (a)).  Since $\E[\kappa(R_n)] <2$ by~(\ref{eqn.kappa0}), 
if $n^*$ is sufficiently large then 
$\pr(\kappa(R_n) \leq n^*) \geq (1+\eps/3)^{-1}$. 
Now if $N^*$ is sufficiently large and $n \geq N^*$, then by summing~(\ref{eqn.probsb}) (with $\eps$ replaced by $\eps/3$) over all $0 \leq n_1,\ldots,n_k \leq n^*$ we obtain
\begin{eqnarray*}
  |\cA_n|
& \leq &
  (1+\eps/3) \,|\{G \in \cA_n : \land_{i=1}^{k} \kappa(G,H_i) \leq n^*\}|\\ 
& \leq &
  (1+\eps/3)^2\, |\cA_n({\bf 0})|\, e^{\sigma_k} \,\prod_{i=1}^k \pr(X_i \leq n^*)\\
& \leq &
  (1+\eps/3)^2\, |\cA_n({\bf 0})|\, e^{\sigma_k} \;\; \leq \; (1+\eps)\, |\cA_n({\bf 0})|\, e^{\sigma_k}\,.
\end{eqnarray*}
Thus
\[ |\cA_n({\bf 0})|/ |\cA_n| \geq (1+\eps)^{-1} e^{-\sigma_k} \geq (1-\eps)\, e^{-\sigma_k}\,. \]
This proves~(\ref{eqn.probsb0}) and thus completes the proof of part~(b).
\end{proof}

The following lemmas (like the last lemma) involve a class $\cA$ of graphs and a class $\cC$ of connected graphs satisfying several conditions, which say roughly that we stay in $\cA$ if we attach or detach graphs in~$\cC$ or bridges incident with induced subgraphs in $\cC$.
The next lemma, Lemma~\ref{lem.Hinew}, will be proved quickly from Lemma~\ref{lem.abnew}. In Lemma~\ref{lem.Hinew}, the relevant integer $k$ is required to be at least some value $k_0$ (depending on $\eps$).  In Lemma~\ref{lem.Hi2new} we drop this restriction on $k$ (using Lemma~\ref{lem.Hinew} in the proof).
For $x,y >0$ and $0<\eps<1$ the notation $y=(1 \pm \eps)x$ means $(1-\eps)x \leq y \leq (1+\eps)x$. 

\begin{lemma} \label{lem.Hinew} 
Let the bridge-addable class 
$\cA$ of graphs have growth constant $\gamma$, and let $\rho =\rho_\cA\, (=1/\gamma)$. Let $\cC$ be a non-empty class of connected graphs, and let each graph in $\cC$ be bridge-deletable in $\cA$. List the (unlabelled) graphs in $\tC$ as $H_1,H_2,\ldots$ so that $v(H_1) \leq v(H_2) \leq \cdots$.
As before, call $k \in \N$ \emph{relevant} when $k \leq |\tC|$ (where $\tC$ may be infinite). For each relevant $i$ let $h_i=v(H_i)$, let $\alpha_i= \rho^{h_i}/ \aut H_i$ and let $X_i \sim \Po(\alpha_i)$; and let $X_1, X_2,\ldots$ be independent.  Also, for each relevant $k$ let $\sigma_k = \sum_{i=1}^k \alpha_i$.
Suppose that every graph $H_i$ in $\tC$ is weakly attachable to $\cA$ and weakly detachable from $\cA$.
Let $R_n \inu \cA$.
 
Then for every $\eps>0$ there exists a relevant $k_0 \in \N$ such that for all relevant $k \geq k_0$ and $n^* \in \N$ the following holds: there exists $N^* \in \N$ such that for all $n \geq N^*$ and all $\,0 \leq n_1,\ldots,n_k \leq n^*$ 
\begin{equation} \label{eqn.probs} 
\pr\left( \land_{i=1}^k \, \kappa(R_n,H_i) = n_i \right) = (1 \pm \eps)\,
\prod_{i =1}^{k} \pr(X_i=n_i) \,.
\end{equation}
\end{lemma}
\begin{proof} 
Let $0<\eps<1$. We will use Lemma~\ref{lem.abnew} with $\eps$ replaced by $\eta = \eps/3$.
Let $k_0 \in \N$  be as given in Lemma~\ref{lem.abnew} (a) (with $\eps$ replaced by $\eta$), and let $k \geq k_0$ be relevant. Let $n^* \in \N$. By Lemma~\ref{lem.abnew} there exists $N^*$ sufficiently large that all of~(\ref{eqn.probsa0}),~(\ref{eqn.probsa}),~(\ref{eqn.probsb0}) and~(\ref{eqn.probsb}) hold for all $n \geq N^*$.  Let $n \geq N^*$ and let $\,0 \leq n_1,\ldots,n_k \leq n^*$.  Then by~(\ref{eqn.probsa}) and~(\ref{eqn.probsb0}), since $(1-\eta)^2 \geq 1-\eps$,
\[ \pr\left( \land_{i=1}^k \, \kappa(R_n,H_i) = n_i \right) \geq (1-\eps) \, \prod_{i =1}^{k} \pr(X_i=n_i)\,; \]
and by~(\ref{eqn.probsb}) and~(\ref{eqn.probsa0}), since $(1+\eta)^2 \leq 1+\eps$,
\[ \pr\left( \land_{i=1}^k \, \kappa(R_n,H_i) = n_i \right) \leq (1 + \eps) \, \prod_{i =1}^{k} \pr(X_i=n_i). \]
Thus~(\ref{eqn.probs}) holds, as required.
\end{proof}

We now show that in the last lemma we do not need to insist that $k$ be large.  The following lemma has the same premises as Lemma~\ref{lem.Hinew}.

\begin{lemma} \label{lem.Hi2new}
  Let the premises be as in Lemma~\ref{lem.Hinew}. Then for every $\eps>0$, relevant $k \in \N$ and $n^* \in \N$ the following holds: there exists $N^* \in \N$ such that for all $n \geq N^*$ and all $\,0 \leq n_1,\ldots,n_k \leq n^*$ 
\begin{equation} \label{eqn.probs2}  
\pr\left( \land_{i=1}^k \, \kappa(R_n,H_i) = n_i \right) = (1 \pm \eps)\, \prod_{i=1}^k \pr(X_i=n_i) \,.
\end{equation}
\end{lemma}
\begin{proof}
Let $\eps>0$, let $k \in \N$ be relevant and let $n^* \in \N$.  We shall use Lemma~\ref{lem.Hinew} with $\eps$ replaced by $\eps/2$ to prove~(\ref{eqn.probs2}).  Let $k_0$ be as given in Lemma~\ref{lem.Hinew} (with $\eps$ replaced by $\eps/2$).  If $k \geq k_0$ we may take $N^*$ as in Lemma~\ref{lem.Hinew} and we are done; 
so we may assume that $k<k_0$. Let
\begin{equation} \label{eqn.mu}
\mu = \min_{n_i \leq n^*} \, \prod_{i=1}^{k} \pr(X_i = n_i)
\end{equation}
where the minimum is over all $0 \leq n_1,\ldots,n_k \leq n^*$. Recall that by~(\ref{eqn.kappa0}) we have $\kappa(R_n) \leq_s 1+\Po(1)$. Let $\hat{n}$ be sufficiently large that
(i) $\pr(\kappa(R_n) > \hat{n}) \leq (\eps/2) \, \mu$ for each $n \in \N$,
and (ii) $\prod_{i=k+1}^{k_0} \pr(X_i \leq \hat{n}) \geq (1-\eps/2)$.  Let $N^*$ be as given by Lemma~\ref{lem.Hinew}, with $\eps$ replaced by 
$\eps/2$, with $k=k_0$ and with $n^*$ replaced by $\hat{n}$.

Let $n \geq N^*$, and let $0\leq n_1,\ldots,n_k \leq n^*$.  By summing in~(\ref{eqn.probs}) over all $0 \leq n_{k+1},\ldots,n_{k_0} \leq \hat{n}$ we obtain
\begin{eqnarray*}
&&
\pr\left( (\land_{i=1}^k \, \kappa(R_n,H_i) = n_i) \land ( \land_{i=k+1}^{k_0} \kappa(R_n,H_i) \leq \hat{n}) \right)\\
&=&
(1 \pm \eps/2) \,
\prod_{i=1}^{k} \pr(X_i=n_i) \; \prod_{i=k+1}^{k_0} \pr(X_i \leq \hat{n}). 
\end{eqnarray*} 
Hence
\begin{eqnarray*}
\pr\left( \land_{i=1}^k \, \kappa(R_n,H_i) = n_i \right)
& \geq &
(1 - \eps/2)^2 \,
\prod_{i=1}^{k} \pr(X_i=n_i)\\
 & \geq &
 (1 - \eps) \,
\prod_{i=1}^{k} \pr(X_i=n_i);
\end{eqnarray*} 
and
\begin{eqnarray*}
&&
\pr\left( \land_{i=1}^k \, \kappa(R_n,H_i) = n_i \right)\\
& \leq &
(1 + \eps/2) \,
\prod_{i=1}^{k} \pr(X_i=n_i) + \pr(\kappa(R_n) >\hat{n})\\
& \leq &
(1 + \eps) \,
\prod_{i=1}^{k} \pr(X_i=n_i). 
\end{eqnarray*} 
Thus~(\ref{eqn.probs2}) holds, as required.
\end{proof}

We need to prove that $C(\rho)$ is finite. The following lemma shows that $C(\rho) \leq 1$, and we shall soon improve the bound to $C(\rho) \leq \frac12$.

\begin{lemma} \label{lem.Crhonew}
Let the premises be as in Lemmas~\ref{lem.Hinew} and~\ref{lem.Hi2new}, and let $C(x)$ be the exponential generating function for $\cC$.  Then $C(\rho) \leq 1$.
\end{lemma}
\begin{proof}
Fix a relevant $k \in \N$, and let $X_n = \sum_{i=1}^k \kappa(R_n,H_i)$ for $n \in \N$. Let $Y \sim \Po(\sigma_k)$, where $\sigma_k = \sum_{i=1}^k \alpha_i$.
Let $\eps>0$.  If $n^*$ is sufficiently large then
\[ \pr(X_n > n^*) \leq \pr(\kappa(R_n) > n^*) < \eps/3\ \;\;\;\; \mbox{ using~(\ref{eqn.kappa0})},\]
and $\pr(Y>n^*) < \eps/3$.  By Lemma~\ref{lem.Hi2new} there exists $N$ such that for all $n \geq N$, for all $A \subseteq \{0,1,\ldots,n^*\}$
\[ |\pr(X_n \in A) - \pr(Y \in A)| < \eps/3. \]
Hence $d_{TV}(X_n,Y) < \eps$ for all $n \geq N$.  Thus  
$X_n$ converges in distribution to $Y$ as $n \to \infty$, and so
\[ \liminf_{n \to \infty}\, \E[X_n] \geq \E[Y] = \sigma_k\,.\] 
%
Also, using~(\ref{eqn.kappa0}) and then~(\ref{eqn.frag})
\begin{eqnarray*}
 2 & \geq & \E[\kappa(R_n)] \geq \E[X_n] + \pr({\rm Big}(R_n) \not\in \{H_1,\ldots,H_k\})\\ 
 & \geq & \E[X_n] + 1 + o(1).
\end{eqnarray*}
It follows that $\sigma_k \leq 1$ for each relevant $k$. But
\[ \sigma_k = \sum_{i=1}^k \alpha_i = \sum_{i=1}^k \frac{\rho^{h_i}}{\aut(H_i)} = \sum_{i=1}^k  \frac{h_i !}{\aut(H_i)}\, \frac{\rho^{h_i}}{h_i !} = \sum_{i=1}^k  c_i \, \frac{\rho^{h_i}}{h_i !}\]
where $c_i$ is the number of graphs isomorphic to $H_i$ on vertex set $[h_i]$. Hence, if $\tilde{\cC}$ is finite and $k=|\tilde{\cC}|$ then $\sigma_k= C(\rho)$; and
if $\tilde{\cC}$ is infinite and $k \to \infty$ then $\sigma_k \to C(\rho)$. Thus $C(\rho) \leq 1$, as required.
\end{proof}

From Lemmas~\ref{lem.Hi2new} and~\ref{lem.Crhonew} 
we obtain the following lemma, Lemma~\ref{lem.prefree}.  Note that $C(\rho)< \infty$ by Lemma~\ref{lem.Crhonew}, so $\BP(\cC,\rho)$ is well defined in Lemma~\ref{lem.prefree}.

\begin{lemma} \label{lem.prefree} 
Let the premises be as in Lemmas~\ref{lem.Hinew} to~\ref{lem.Crhonew}. Then $\Frag(R_n,\cC)$ converges in total variation to $\BP(\cC,\rho)$ as $n \to \infty$.
\end{lemma}
\begin{proof}
 Write $F_n$ for $\Frag(R_n,\cC)$ for $n \in \N$, and
let $R \sim \BP(\cC,\rho)$.  Let $\cF$ be the class of graphs such that each component is in $\cC$. Observe that the set $\tF$ of unlabelled graphs is countable, and our random variables take values in this set. 
Let $\eps>0$.  We must show that there exists $N^*$ such that for every $A \subseteq \tF$ and every $n \geq N^*$ we have
\begin{equation} \label{eqn.totvar1}
\big| \pr(F_n \in A) - \pr(R \in A) \big| < \eps.
\end{equation}
Since $v(F_n) \leq \frag(R_n)$ and $\E[\frag(R_n)]<2$ by~(\ref{eqn.frag}), there exists $n^*$ such that $\pr(v(F_n)> n^*) < \eps/2$ and $\pr(v(R)> n^*) < \eps/2$.  Let $F_n^-$, $R^-$ be the subgraph of $F_n$, $R$ respectively consisting of the components of order at most~$n^*$.  If $A^-$ is the set of graphs in $A$ of order at most $n^*$ then
\[ \big| \pr(F_n \in A) - \pr(R \in A) \big| < 
\big| \pr(F^-_n \in A^-) - \pr(R^- \in A^-) \big| + \eps/2\,. \]
Thus to prove~(\ref{eqn.totvar1}) it suffices to show that there exists $N^*$ such that for every non-empty subset $A$ of graphs in $\tF$ each of order at most $n^*$ and every $n \geq N^*$ we have
\begin{equation} \label{eqn.totvar2} 
\big| \pr(F_n^- \in A) - \pr(R^- \in A) \big| \leq \eps/2.
\end{equation}

List the (connected, unlabelled) graphs in $\tC$ as $H_1,H_2,\ldots$ such that $v(H_1) \leq v(H_2) \leq \cdots$.  If $\tilde{\cC}$ is finite let $k=|\tilde{\cC}|$; and otherwise let $k$ be the least positive integer such that $v(H_{k+1})> n^*$.
Let $\mu$ be as in equation~(\ref{eqn.mu}).  
Since $\cA$ has a growth constant, and $n^*$ is fixed, there exists $N^*_0$ such that 
\[ \pr(\giant(R_n) \leq n^*)  \leq \eps \mu /5 \;\; \mbox{ for each } n \geq N^*_0\, .\] 
(We do not need here to use the assumption that $\cA$ is bridge-addable.) 
By Lemma~\ref{lem.Hi2new} (with $\eps$ replaced by $\eps/5$) there exists $N^* \geq N^*_0$ such that for all $n \geq N^*$ and all $\,0 \leq n_1,\ldots,n_k \leq n^*$  
\begin{equation} \label{eqn.kappa2}
\pr\left( \land_{i=1}^k \, \kappa(R_n,H_i) = n_i \right) = (1 \pm \eps/5)\,
\prod_{i =1}^{k} \pr(X_i=n_i) \,.
\end{equation}
Each graph $G \in A$ corresponds to a unique $k$-tuple $\bn = (n_1,\ldots,n_k)$ with $\,0 \leq n_1,\ldots,n_k \leq n^*$.  Let $B$ be the (non-empty) set of $k$-tuples $\bf n$ corresponding to the graphs in $A$.
Let $n \geq N^*$. Then
\[ \pr(F_n^- \in A) = \sum_{\bn \in B} \pr\big(\land_{i=1}^{k} \kappa(\Frag(R_n),H_i) = n_i)\,,\]
so
\[ \big| \pr(F_n^- \in A) - \sum_{\bn \in B} \pr\big(\land_{i=1}^{k} \kappa(R_n,H_i) =n_i \big)\big| \leq \pr(\giant(R_n) \leq n^*) \leq \eps \mu /5\,,\]
and thus
\[ \pr(F_n^- \in A) =  (1 \pm \eps/5) \, \sum_{\bn \in B} \pr\big(\land_{i=1}^{k} \kappa(R_n,H_i) =n_i \big) \,.\]
Hence by~(\ref{eqn.kappa2})
\begin{eqnarray*}
  \pr(F_n^- \in A)
&=&
  (1 \pm \eps/5)^2 \, \sum_{\bn \in B}\, \prod_{i=1}^{k} \, \pr(X_i=n_i)\\
  &=&
  (1 \pm \eps/2) \; \pr(R^- \in A)\,.
\end{eqnarray*}
The inequality~(\ref{eqn.totvar2}) follows, and this completes the proof.
\end{proof}

Assume that the conditions in Lemmas~\ref{lem.Hinew} to~\ref{lem.prefree} hold. Then by Lemma~\ref{lem.prefree}, 
\[ \pr(\Frag(R_n,\cC)=\emptyset) \to e^{-C(\rho)}\]
as $n \to \infty$.  By~(\ref{eqn.badd-conn12}) (quoted from reference~\cite{Chapuy-Perarnau-2019})
\[ \pr(\Frag(R_n,\cC)=\emptyset) \geq \pr(R_n \mbox{ is connected}) \geq e^{-\frac12} +o(1)\,.\]
Hence $C(\rho) \leq \frac12$, improving on Lemma~\ref{lem.Crhonew}. We have now completed the proof of Theorem~\ref{thm.comps-badd}.  
\medskip

Next let us prove the result~(\ref{eqn.Cdash}) on 
$\E[v(R)]$ stated after Theorem~\ref{thm.comps-badd}.
\begin{proof}[Proof of~(\ref{eqn.Cdash})]
By~(\ref{eqn.frag}) and Theorem~\ref{thm.comps-badd}, for each $k \in \N$
\[ 2 > 
\sum_{i=1}^{k} i\, \pr\big(v(\Frag(R_n,\cC))=i\big) \rightarrow 
\sum_{i=1}^{k} i\, \pr(v(R)=i) \;\;\; \mbox{ as } n \to \infty\,.\]
Thus $\sum_{i=1}^{k} i\,\pr(v(R)=i) \leq 2$, and so $\E[v(R)] \leq 2$.  Also, see~(\ref{eqn.orderk}),
\[ \E[v(R)] = \sum_{k \geq 1} |\cF_k| \frac{k \rho^{k-1}}{k!} \cdot \frac{\rho}{F(\rho)} = F'(\rho) \cdot \frac{\rho}{F(\rho)} = \rho\, C'(\rho)\,.\]
This completes the proof of~(\ref{eqn.Cdash}).
\end{proof}
\smallskip

Finally in this section we supply the missing steps in our proof of Corollary~\ref{cor.cycle}, first concerning the growth constant and then concerning $\Frag(R_n)$ being planar.
\begin{proof}[Proof that, when $\cA$ is $\cE^S$, $\cA^t$ in Corollary~\ref{cor.cycle} has growth constant $\gamma_\cP$]
We may assume that $t$ is $o(n)$. Let $\cC$ be the class of connected planar graphs. For simplicity assume that $t | n$, and let $k=n/t$.  We construct many connected graphs in $\cA_n^{\,t}$ as follows.
Partition $[n]$ into $k$ $t$-sets $V_1,\ldots,V_k$.  Let $x_i \in V_i$ for $i=1,\ldots,k$; and insert edges $x_ix_{i+1}$ for $i=1,\ldots,k\!-\!1$ to form a path.
Finally put a $t$-vertex graph in $\cC$ on each set $V_i$.  This forms a connected planar graph in $\cA^t_n$.  The number of constructions is
\[ \frac{n!}{(t!)^k}\, t^k\, |\cC_t|^k 
= n! \, \left(|\cC_t| / (t\!-\!1)!\right)^{n/t}.\]

Each of these graphs $G$ is constructed just twice (assuming $k \geq 2$).  For, in $G$ there are exactly $k-1$ bridges $b$ such that when we delete $b$ both components have order divisible by $t$, and these bridges form a path.  Given $G$ and the correct orientation of the path we can read off the entire construction.  Hence the number of distinct graphs constructed is 
$\tfrac12 \, n! \, ( |\cC_t| / (t\!-\!1)!)^{n/t}$.
But $|\cC_t|/ (t\!-\!1)! = (\gamma_\cP +o(1))^t$ as $n \to \infty$ (and so also $t \to \infty$), so
\[  |\cA^{\,t}_n| \geq n! \, (\gamma_\cP+o(1))^n\,.\]
Since $\cA^t \subseteq \cA$ it follows 
that $\cA^t$ has growth constant $\gamma_\cP$, as required.
\end{proof}

\begin{proof}[Proof that, when $\cA$ is $\cE^S$, whp $\Frag(R_n)$ is planar]

Let $0 < \eps < 1$, and let $c = 4/\eps$.
Then $\pr(\frag(R_n) \geq c) \leq \eps/2$ by~(\ref{eqn.frag}).
Let $\cH$ be the (finite) set of unlabelled non-planar connected graphs in $\cE^S$ of order less than~$c$. For each $H \in \cH$ let $\cB^H$ be the class of graphs $G$ in $\cA^{t}$ such that $\frag(G) < c$ and $\Frag(G)$ has a component 
$H$, and let $\cB$ be the union of these classes $\cB^H$ over $H \in \cH$.
Suppose that $n \geq 2c$. 

From each graph $G$ in $\cB^H_n$, by adding an edge between a component $H$ and a vertex in Big$(G)$ we construct at least $v(H) (n - c) \geq n/2$ graphs $G'$ in $\cA^{t}_n$.  Note that big$(G') > c + v(H) > 2\, v(H)$, so any pendant appearances of $H$ in Big$(G')$ are vertex-disjoint.
Let $g$ be the Euler genus of the surface $S$. Then each constructed graph $G'$ has at most $g$ pendant appearances of $H$ in its big component, so $G'$ can be constructed at most $g$ times.  Hence
 $|\cB^H_n| \, n/2 \leq |\cA^{t}_n| \, g$,
and so
$ |\cB_n| \, n/2 \leq |\cA^{t}_n| \, g \cdot |\cH|$.  Thus
\[ \pr(R_n \in \cB) = |\cB_n| / |\cA^{t}_n| \leq 2 g \,|\cH|/n =o(1).\]
Finally,
\[ \pr(\Frag(R_n) \mbox{ is non-planar}) \leq \pr(\frag(R_n) \geq  c) + \pr(R_n \in \cB) \leq \eps/2+o(1);\]
and so $\Frag(R_n)$ is planar whp, as required.
\end{proof}

\section{Concluding Remarks}
\label{sec.concl}

We investigated random graphs $R_n$ sampled uniformly from a structured class $\cA$ of graphs, such as the class $\cE^S$ of graphs embeddable in a given surface~$S$ or an addable minor-closed class. We sharpened results on pendant appearances in $R_n$ when the class has a growth constant (or at least a finite non-zero radius of convergence) to give the correct constant.

When $\cA$ has a growth constant and is bridge-addable and $\cC$ satisfies certain conditions, we deduced results concerning the fragment $\Frag(R_n,\cC)$, see Theorem~\ref{thm.comps-badd} and its Corollaries, in particular determining the limiting distribution of $\Frag(R_n)$ in certain cases and thus the limiting probability of $R_n$ being connected.
Corollary~\ref{cor.cycle}
has an additional constraint, namely an upper bound on cycle lengths: 
with further work~\cite{cmcd-badd-Delta}, we can obtain a similar result 
when the additional constraint is an upper bound on vertex degrees. 
%
As was noted earlier (in Section~\ref{subsec.relwork}), related results under stronger assumptions can be obtained by first proving that the relevant class is smooth, and in this case we may also learn for example about the core of $R_n$, see~\cite{cmcd-rgstrcl}.

Finally, note that throughout we have been concerned with random graphs $R_n$ from a class of \emph{labelled} graphs.  We know very little for unlabelled graphs.
For example, let $\cA$ be a bridge-addable class of labelled graphs.  It was useful above that for  $R_n \in_u \cA$ we have inequalities like~(\ref{eqn.kappa1}) and~(\ref{eqn.badd-conn12}).  But let $\tA$ be the set of unlabelled graphs in $\cA$, and let $\tilde{R}_n$ be uniformly distributed over the $n$-vertex graphs in $\tA$.  Is there a $\delta>0$ such that $\pr(\tilde{R}_n \mbox{ is connected}) \geq \delta$ for all~$n$?  See~\cite{cmcd-unlab2020} concerning this and related conjectures.
\bigskip

\noindent{\bf Acknowledgement} \;
I would like to thank the referees for positive and helpful comments.



\begin{thebibliography}{99}


\bibitem{bcr08}
E.A.~Bender, E.R.~Canfield and L.B.~Richmond.
Coefficients of functional compositions often grow smoothly,
{\em Electron. J. Combin.} {\bf 15} (2008) \#R21.


\bibitem{BNW} O.~Bernardi, M.~Noy and D.J.A.~Welsh,
On the growth rate of minor-closed classes of graphs,
\textit{J. Combin. Theory Ser. B} {\bf 100} (2010), 468 -- 484.


\bibitem{BM-KW-2014}
M.~Bousquet-M\'elou and K.~Weller.
Asymptotic Properties of Some Minor-Closed Classes of Graphs,
\textit{Combinatorics Probability and Computing} {\bf 23} (2014), 749 -- 795.   

\bibitem{CFGMN2011}
G.~Chapuy, E.~Fusy, O.~Gim\'{e}nez, B.~Mohar, and M.~Noy. 
Asymptotic enumeration and limit laws for graphs of fixed genus, 
\textit{J. Combin. Theory Ser. A} {\bf 118} (2011), 748--777.

\bibitem{Chapuy-Perarnau-2019} G.~Chapuy and G.~Perarnau. 
Connectivity in bridge-addable graph classes: the McDiarmid-Steger-Welsh conjecture, \emph{J. Comb. Th. B} {\bf 136} (2019), 44--71.

\bibitem{Chapuy-Perarnau-2020} G.~Chapuy and G.~Perarnau. 
Local convergence and stability of tight bridge-addable graph classes, \emph{Canadian Journal of Mathematics} {\bf 72} (2020), 563 -- 601.


\bibitem{Diestel}
R.~Diestel.  \emph{Graph Theory}, 5th edition, Springer, 2017.




\bibitem{GN2009}
O.~Gim\'{e}nez and M.~Noy. 
Asymptotic enumeration and limit laws of planar graphs, 
\emph{J.~Amer. Math. Soc.} {\bf 22} (2009), 309--329.

\bibitem{KP13}
M.~Kang and K.~Panagiotou.
On the connectivity of random graphs
from addable classes, \emph{J. Combin. Theory Ser. B}, {\bf 103}(2) (2013), 306--312.

\bibitem{msw2005} C.~McDiarmid, A.~Steger and D.J.A.~Welsh.
Random planar graphs,
\textit{J.~Combin. Theory Ser. B} \textbf{93} (2005),  187--205.

\bibitem{msw2006} C.~McDiarmid, A.~Steger and  D.J.A.~Welsh. 
Random graphs from planar and other addable classes. \emph{Topics in Discrete Mathematics,} 231--246, Algorithms Combin. 26, Springer, Berlin, 2006.

\bibitem{cmcd-rgons}
C.~McDiarmid.
Random graphs on surfaces, \textit{J. Combin. Theory Ser. B} {\bf 98} (2008), 778 -- 797.

\bibitem{cmcd-rgmc}
C.~McDiarmid.
Random graphs from a minor-closed class,
{\em Combinatorics, Probability and Computing} {\bf 18} (2009), 583-599.

\bibitem{cmcd-connwba2012}
C.~McDiarmid.
Connectivity for random graphs from a weighted bridge-addable class,
{\em Electronic J. Combinatorics} {\bf 19}(4) (2012), P53.

\bibitem{cmcd-rgwmc2013}
C.~McDiarmid.
Random graphs from a weighted minor-closed class,
\emph{Electronic J. Combinatorics} {\bf20} (2) (2013) P52, 39 pages.
  

\bibitem{cmcd-bham}
C.~McDiarmid.
On random graphs from a minor-closed class,
Chapter 4, pages 102 -- 120, in: \emph{Random Graphs, Geometry and Asymptotic Structure}, LMS Student Texts 84, 2016.


\bibitem{cmcd-PAC}
C.~McDiarmid. 
Pendant appearances and components in random graphs from structured classes, arXiv:2108.07051, 2021.

\bibitem{cmcd-unlab2020} C.~McDiarmid. Connectivity for an unlabelled bridge-addable graph class,
Open problem, First Armenian Workshop on Graphs, Combinatorics, Probability, May 2019.
arXiv:2001.05256v2. 

\bibitem{cmcd-rgstrcl}
C.~McDiarmid.
Random graphs from structured classes,
arXiv: 2209.10476v1, 2022.


\bibitem{cmcd-badd-Delta} C.~McDiarmid. Random graphs, bridge-addability and bounded degrees, in preparation, 2023.


\bibitem{ColinSophiaSizes}
C.~McDiarmid and S. Saller. Classes of graphs embeddable in order-dependent surfaces, \textit{Combinatorial Theory} 3.1 (2023).


\bibitem{ColinSophiaProperties}
C.~McDiarmid and S. Saller. Random graphs embeddable in order-dependent surfaces, 
\emph{Random Structures and Algorithms} (2023), 1 -- 46. (online 5/12/2023).




\bibitem{mt2001}
B.~Mohar and C.~Thomassen, 
\emph{Graphs on Surfaces},
Johns Hopkins University Press, Baltimore, MD (2001).

\bibitem{Moon70}
J.W. Moon. \emph{Counting labelled trees}, Vol 1 of \emph{Canadian Mathematical Monographs}, 1970.

\bibitem{Renyi59}
A. R\'{e}nyi. Some remarks on the theory of trees, \emph{Publications of the Mathematical Institute of the  Hungarian Academy of Sciences} {\bf 4} (1959), 73 -- 85. 

\bibitem{Stufler-Gibbs}
B. Stufler.  Gibbs partitions: the convergent case, 
\emph{Random Structures and Algorithms} {\bf 53} (2018), 537 -- 558.

\bibitem{Stufler-survey}
B. Stufler. Limits of random tree-like discrete structures, \emph{Probability Surveys} {\bf 17} (2020), 318–477.

\end{thebibliography}
\end{document}